\documentclass[a4paper,12pt]{amsart}

\usepackage{amsmath}
\usepackage[latin1]{inputenc}
\usepackage{amsfonts}
\usepackage{amssymb}
\usepackage{geometry}
\newtheorem{thm}{Theorem}[section]
\newtheorem{cor}[thm]{Corollary}
\newtheorem{lem}[thm]{Lemma}

\newtheorem{defn}[thm]{Definition}
\newtheorem{rem}[thm]{Remark}

\numberwithin{equation}{section}

\def\R{\mathbb R}

\def\b{\mathrm{b}}
\def\a{a}
\def\B{B}

\newcommand{\mand}{\quad\text{and}\quad}

\title[General Reiteration Theorems for ${\mathcal R}$ and ${\mathcal L}$ classes]
{General Reiteration Theorems for ${\mathcal R}$ and ${\mathcal L}$ classes: Mixed interpolation of ${\mathcal R}$ and ${\mathcal L}$-spaces}

\author[Fern\'andez-Mart\'{\i}nez]{Pedro Fern\'andez-Mart\'{\i}nez}
\address[Pedro Fern\'andez-Mart\'{\i}nez and Teresa M. Signes]{Departamento de Matem\'aticas \\
Facultad de Matemáticas \\ Universidad de Murcia \\ Campus de
Espinardo \\ 30071 Espinardo (Murcia), Spain}

\author[Signes]{Teresa M. Signes}

\date{\today}

\begin{document}
\begin{abstract}
Given $E_0, E_1, F_0, F_1, E$ rearrangement invariant function spaces, $a_0$, $a_1$, $\b_0$, $\b_1$, $\b$ slowly varying functions and $0< \theta_0<\theta_1<1$, we characterize the interpolation spaces
$$(\overline{X}^{\mathcal R}_{\theta_0,\b_0,E_0,\a_0,F_0}, \overline{X}^{\mathcal R}_{\theta_1, \b_1,E_1,\a_1,F_1})_{\theta,\b,E},\quad
(\overline{X}^{\mathcal L}_{\theta_0, \b_0,E_0,\a_0,F_0}, \overline{X}^{\mathcal L}_{\theta_1,\b_1,E_1,\a_1,F_1})_{\theta,\b,E}$$
and
$$(\overline{X}^{\mathcal R}_{\theta_0,\b_0,E_0,\a_0,F_0}, \overline{X}^{\mathcal L}_{\theta_1, \b_1,E_1,\a_1,F_1})_{\theta,\b,E},\quad
(\overline{X}^{\mathcal L}_{\theta_0, \b_0,E_0,\a_0,F_0}, \overline{X}^{\mathcal R}_{\theta_1,\b_1,E_1,\a_1,F_1})_{\theta,\b,E},$$
for all possible values of $\theta\in[0,1]$. Applications to interpolation identities for grand and small Lebesgue spaces, Gamma spaces and $A$ and $B$-type spaces are given.
\end{abstract}
\maketitle

\section{Introduction}\label{introduction}


This paper continues the study initiated in \cite{FMS-RL1,FMS-RL2} where the present authors proved reiteration theorems for couples formed by the spaces
$$\overline{X}_{\theta,\b,E}, 
\quad \overline{X}^{\mathcal R}_{\theta,\b,E,\a,F},\quad \overline{X}^{\mathcal L}_{\theta,\b,E,\a,F}$$
where $\theta\in[0,1]$, $a$ and $b$ are slowly varying functions and $E$, $F$ are rearrangement invariant (r.i.) function  spaces. Given a (quasi-) Banach couple $\overline{X}=(X_0,X_1)$ these spaces are defined as 
$$
 \overline{X}_{\theta,\b,E}=\Big\{f\in X_0+X_1\;\colon\;
\big \| t^{-\theta} {\b}(t) K(t,f) \big \|_{\widetilde{E}} < \infty\Big\},
$$
$$\overline{X}_{\theta,\b,E,a,F}^{\mathcal R}=\Big\{f\in X_0+X_1\;\colon\; \Big \|  \b(t) \|   s^{-\theta} a(s) K(s,f) \|_{\widetilde{F}(t,\infty)}      \Big   \|_{\widetilde{E}}<\infty\Big\}$$
and
$$\overline{X}_{\theta,\b,E,a,F}^{\mathcal L}=\Big\{f\in X_0+X_1\;\colon\;
\Big\|\b(t) \|s^{-\theta} a(s) K(s,f) \|_{\widetilde{F}(0,t)}\Big\|_{\widetilde{E}} < \infty\Big\},$$
where $K(\cdot,f)$ is the Peetre $K$-functional and $\widetilde{E}$ is the corresponding r.i. space with respect the homogeneous measure $dt/t$ on $(0,\infty)$; see \S 2 below.

In \cite{FMS-RL1,FMS-RL2}, we characterized the interpolation spaces 
\begin{equation}\label{ec1}
\big(\overline{X}^{\mathcal R}_{\theta_0, \b_0,E_0,\a,F}, \overline{X}_{\theta_1,\b_1,E_1}\big)_{\theta,\b,E},\qquad
\big(\overline{X}_{\theta_0,\b_0,E_0}, \overline{X}^{\mathcal L}_{\theta_1, \b_1,E_1,\a,F}\big)_{\theta,\b,E}
\end{equation}
and its ``dual'' situation
\begin{equation}\label{ec2}
\big(\overline{X}_{\theta_0,\b_0,E_0}, \overline{X}^{\mathcal R}_{\theta_1, \b_1,E_1,\a,F}\big)_{\theta,\b,E},\qquad
\big(\overline{X}^{\mathcal L}_{\theta_0, \b_0,E_0,\a,F}, \overline{X}_{\theta_1,\b_1,E_1}\big)_{\theta,\b,E}.
\end{equation}
On this occasion we shall focus on the reiteration results for couples formed only by ${\mathcal R}$ and ${\mathcal L}$-spaces, specifically
\begin{equation}\label{ec3}
(\overline{X}^{\mathcal R}_{\theta_0,\b_0,E_0,\a_0,F_0}, \overline{X}^{\mathcal R}_{\theta_1, \b_1,E_1,\a_1,F_1})_{\theta,\b,E},\quad
(\overline{X}^{\mathcal L}_{\theta_0, \b_0,E_0,\a_0,F_0}, \overline{X}^{\mathcal L}_{\theta_1,\b_1,E_1,\a_1,F_1})_{\theta,\b,E}
\end{equation}
and
\begin{equation}\label{ec4}
(\overline{X}^{\mathcal R}_{\theta_0,\b_0,E_0,\a_0,F_0}, \overline{X}^{\mathcal L}_{\theta_1, \b_1,E_1,\a_1,F_1})_{\theta,\b,E},\quad
(\overline{X}^{\mathcal L}_{\theta_0, \b_0,E_0,\a_0,F_0}, \overline{X}^{\mathcal R}_{\theta_1,\b_1,E_1,\a_1,F_1})_{\theta,\b,E},
\end{equation}
for all values of $\theta\in[0,1]$.  
We shall show that, when $\theta\in(0,1)$, all these spaces can be described within 
 the classical scale $\overline{X}_{\theta,\b,E}$. However, in the cases $\theta=0, 1$, the resulting reiteration spaces belong to the extremal scales
\begin{equation}\label{extreme1}
\overline{X}^{\mathcal R,\mathcal R}_{\theta,c,E,\b,F,\a,G},\quad \overline{X}^{\mathcal L,\mathcal L}_{\theta,c,E,\b,F,\a,G},\quad \overline{X}^{\mathcal R,\mathcal L}_{\theta,c,E,\b,F,\a,G}\mand  \overline{X}^{\mathcal L,\mathcal R}_{\theta,c,E,\b,F,\a,G},
\end{equation}
introduced in \cite{FMS-RL1,FMS-RL2}; see Definition \ref{defLRR} below. The precise statements are given in Theorems \ref{thmRR} through \ref{thmLR} in \S\ref{reiterationtheorems}.

\

In the special case that $\theta\in(0,1)$ and $E=L_q$, a recent work of Doktorskii \cite{Do2020} provides identities for all the spaces in \eqref{ec3} and \eqref{ec4} (even in the extended quasi-Banach range $0<q\leq \infty$), using a direct and elegant reiteration method. Our approach
in this paper is slightly different, and has the advantage of being applicable also in the extreme cases $\theta=0,1$. Additionally it provides explicit identities for the $K$-functionals of all cases in the form of
Holmsted type formulae that have an independent interest; see \S \ref{sectionH} below.

\


As in \cite{FMS-RL1,FMS-RL2}, the present work finds its motivation in several recent applications to the interpolation of 
 \textit{grand} and \textit{small Lebesgue} spaces $L^{p),\alpha}$, $L^{(p,\alpha}$, $1<p<\infty$, $\alpha>0$ (see Definition \ref{def_gLp} below and the references \cite{AFH,AFFGR,FFGKR}). As  observed in \cite{FK,op1}, these  can be identified as ${\mathcal R}$ and ${\mathcal L}$-spaces in the following way
$$
L^{p),\alpha}=(L_1,L_\infty)^{\mathcal R}_{1-\frac1p,\ell^{-\frac{\alpha}{p}}(t),L_\infty,1,L_p} \mand
L^{(p,\alpha}=(L_1,L_\infty)^{\mathcal L}_{1-\frac1p,\ell^{-\frac{\alpha}{p}+\alpha-1}(t),L_1,1,L_p}
$$
where $\ell(t)=1+|\log(t)|$, $t\in(0,1)$. Then, our results provide identities for
the interpolation spaces
\begin{equation}\label{ec6}
\big(L^{p_0),\alpha},L^{p_1),\beta}\big)_{\theta,\b,E},\qquad\big(L^{(p_0,\alpha},L^{(p_1,\beta}\big)_{\theta,\b,E},\qquad \big(L^{(p_0,\alpha},L^{p_1),\beta}\big)_{\theta,\b,E}
\end{equation}
and also 
\begin{equation}\label{ec7}
\big(L^{p_0),\alpha},L^{(p_1,\beta}\big)_{\theta,\b,E}
\end{equation}
when $1\leq p_0<p_1\leq \infty$, $\alpha, \beta>0$ and $\theta\in[0,1]$. Some of these cases have been obtained in the recent papers \cite{AFH,AFFGR,Do2020,FFGKR} when $\theta\in(0,1)$ and $E=L_q$, $0<q\leq \infty$. 
However, in the limiting cases $\theta=0,1$ the identities for \eqref{ec7} are new, and those for \eqref{ec6} have been improved (with respect to \cite{FMS-RL1,FMS-RL2}) in the sense that we can now identify the resulting spaces only in terms of the couple $(L_1,L_\infty)$; see  Corollaries \ref{cor55} c), \ref{cor57} b), \ref{cor2} b).


Our formulae also apply to other general $\mathcal{R}$ and $\mathcal{L}$-families such as the Gamma spaces (see \cite{FFGKR}) and the $A$ and $B$-type spaces (see \cite{pu2}). In particular, for these families we can also improve our earlier interpolation identities from \cite[\S5]{FMS-RL2}
 in the same sense as before; see Corollary \ref{corAB} b).

\

The paper is organized as follows. In Section 2 we recall basic concepts regarding rearrangement invariant spaces and slowly varying functions.  We also describe the  interpolation methods we shall work with, namely $\overline{X}_{\theta,\b,E}$, the ${\mathcal R}$ and ${\mathcal L}$-type spaces and the extremal constructions in \eqref{extreme1}.
Generalized Holmstedt type formulae for the $K$-functional  of the different couples  can be found in Section 3. The statements and proofs of reiteration theorems appear in Section 4. Section 5 is devoted to applications.

\section{Preliminaries} \label{preliminaries}

We refer to the monographs \cite{Bennett-Sharpley,Bergh-Lofstrom,B-K,KPS,triebel} for basic concepts on Interpolation Theory and Banach function spaces. A Banach function space $E$ on $(0,\infty)$ is called \textit{rearrangement invariant} (r.i.) if, for any two measurable fuctions $f$, $g$,
\begin{equation*}
g\in E \text{ and } f^*\leq g^*  \Longrightarrow f\in E \text{ and } \|f\|_E\leq \|g\|_E,
\end{equation*}
where $f^*$ and $g^*$ stand for the non-increasing rearrangements of $f$ and $g$.
Following  \cite{Bennett-Sharpley},  we assume that every Banach function space $E$ enjoys the \textit{Fatou property}. Under this assumption  every r.i. space
$E$ is an exact interpolation space with respect to the Banach couple $(L_{1}, L_{\infty})$.

Along this paper we will handle two different measures on $(0,\infty)$: the usual Lebesgue measure $dt$ and the homogeneous measure $\tfrac{dt}{t}$.
We use a tilde to denote r.i. spaces with respect to the second measure.
For example
$$  \|f\|_{\widetilde{L}_{1}} = \int_{0}^{\infty}  |f(t)|  \frac{dt}{t}\quad\mbox{and}\quad \|f\|_{\widetilde{L}_\infty}=\|f\|_{L_\infty}.$$
If $E$ is an r.i. obtained by applying the interpolation functor $\mathcal{F}$ to the couple $(L_1,L_\infty)$, $E=\mathcal{F}(L_1,L_\infty)$, we will denote by $\widetilde{E}$ the space generated by the same $\mathcal{F}$ acting on the couple $(L_1,L_\infty)$, $\widetilde{E}=\mathcal{F}(\widetilde{L}_{1}, L_{\infty})$. 

In the rest of the paper we will denote by $\overline{f}$ the measurable function defined by $$\overline{f}(t)=f(1/t),\quad t>0.$$ Moreover, given $t>0$, it is sometimes convenient to divide the interval $(0,\infty)$ into two subintervals $(0,t)$ and $(t,\infty)$, considering the spaces $\widetilde{E}(0,t)$ and $\widetilde{E}(t,\infty)$ separately. Observe that $f\in \widetilde{E}(t,\infty)$ if and only if $\overline{f}\in \widetilde{E}(0,1/t)$ and
\begin{equation}\label{fbar}
\|f\|_{\widetilde{E}(t,\infty)}=\|\overline{f}\|_{\widetilde{E}(0,1/t)}. 
\end{equation}

Throughout the paper, given two (quasi-) Banach spaces $X$ and $Y$, we will write $X=Y$ if $X\hookrightarrow Y$ and $Y\hookrightarrow X$, where the latter means that $Y\subset X$ and the natural embedding is continuous.
Similarly, $f\sim g$ means that  $f\lesssim g$ and $g\lesssim f$, where $f\lesssim g$ is the abbreviation of $f(t) \leq C g(t)$, $t >0$, for some positive constant $C$ independent of $f$ and $g$.

\subsection{Slowly varying functions}\hspace{1mm}

\vspace{2mm}
In this subsection we recall the definition and basic properties of \textit{slowly varying functions}.  See  \cite{BGT,Matu}.

\begin{defn}\label{def1}
A positive Lebesgue measurable function $\b$, $0\not\equiv\b\not\equiv\infty$,
is said to be \textit{slowly varying} on $(0,\infty)$ (notation $b\in SV$) if, for each $\varepsilon>0$, the function $t \leadsto t^\varepsilon\b(t)$  is equi\-va\-lent to a non-decreasing function on $(0,\infty)$
and $ t \leadsto t^{-\varepsilon}\b(t)$ is equivalent to a non-increasing function on $(0,\infty)$.
\end{defn}

Examples of $SV$-functions include powers of logarithms,
$$\ell^\alpha(t)=(1+|\log t|)^\alpha(t),\quad t>0,\quad \alpha\in\R,$$ ``broken" logarithmic functions defined as
\begin{equation}\label{ellA}
\ell^{(\alpha, \beta)}(t)=\left\{\begin{array}{ll}
\ell^\alpha(t),& 0<t\leq1 \\
\ell^\beta(t)
,& t>1
\end{array}\right.,\quad (\alpha, \beta)\in\R^2,
\end{equation}
reiterated logarithms $(\ell\circ\ldots\circ\ell)^\alpha,\ \alpha\in\R,\ t>0$ and also functions as $t\leadsto\exp(|\log t|^\alpha)$, $\alpha \in (0,1)$. 

Some basic properties of slowly varying functions are summarized in the following lemmas.
\begin{lem}\label{lem0}
Let $\b, \b_1, \b_2\in SV$,  and let $E$ be an r.i. space.
\begin{itemize}
\item[(i)]  Then $\b_1\b_2\in SV$, $\overline{\b}\in SV$ and $\b^r\in SV$ for all $r\in\R$.
\item[(ii)] If $\alpha>0$, then $\b(t^\alpha \b_1(t))\in SV$.
\item[(iii)] If $\epsilon,s>0$ then there are positive constants $c_\epsilon$ and $C_\epsilon$ such that
$$c_\epsilon\min\{s^{-\epsilon},s^\epsilon\}\b(t)\leq \b(st)\leq C_\epsilon \max\{s^\epsilon,s^{-\epsilon}\}\b(t)\quad \mbox{for every}\quad t>0.$$
\end{itemize}
\end{lem}

\begin{lem} \label{lem1}
Let $E$ be an r.i. space and $\b\in SV$.
\begin{itemize}
\item[(i)] If $\alpha>0$, then, for all $t>0$,
$$\|s^\alpha\b(s)\|_{\widetilde{E}(0,t)}\sim t^\alpha\b(t)\quad \mbox{ and } \quad
\|s^{-\alpha}\b(s)\|_{\widetilde{E}(t,\infty)}\sim t^{-\alpha}\b(t).$$
\item[(ii)] If $\alpha\in\R$, then, for all $t>0$,
$$\|s^\alpha\b(s)\|_{\widetilde{E}(t,2t)}\sim t^\alpha\b(t).$$
\item[(iii)] The following functions are slowly varying
$$\Phi_0(t):=\|\b(s)\|_{\widetilde{E}(0,t)}
\quad \text{and}\quad
\Phi_\infty(t):=\|\b(s)\|_{\widetilde{E}(t,\infty)}, \quad t >0.$$
\item[(iv)] For all $t >0$,
$$\b(t) \lesssim \| \b(s) \|_{\widetilde{E}(0,t)} \quad \text{and} \quad
\b(t) \lesssim \| \b(s) \|_{\widetilde{E}(t,\infty)}.$$
\end{itemize}
\end{lem}

We refer to \cite{GOT,FMS-1} for the proof of Lemma \ref{lem0} and \ref{lem1}, respectively.
\begin{rem}\label{rem23}
\textup{The property (iii) of Lemma \ref{lem0} implies that if $\b\in SV$ is such that $\b(t_0)=0$ ($\b(t_0)=\infty$) for some $t_0>0$, then $\b\equiv0$ ($\b\equiv\infty$). Thus, by Lemma \ref{lem1} (iii), if $\|\b\|_{\widetilde{E}(0,1)}<\infty$ then
$\|\b\|_{\widetilde{E}(0,t)}<\infty$ for all $t>0$ and if $\|\b\|_{\widetilde{E}(1,\infty)}<\infty$ then
$\|\b\|_{\widetilde{E}(t,\infty)}<\infty$ for all $t>0$.}

\textup{Moreover, if $f\sim g$ then, using Definition \ref{def1} and Lemma \ref{lem0} (iii), one can show that  $\b\circ f\sim \b\circ g$ for any $b\in SV$.}
\end{rem}

\subsection{Interpolation Methods} \label{interpolation methods}\hspace{2mm}

\vspace{2mm}
Everywhere below $\overline{X}=(X_{0}, X_{1})$ is a \textit{compatible (quasi-) Banach couple}, that is, two  (quasi-) Banach spaces continuously embedded in a Hausdorff topological vector space. The \textit{Peetre $K$-funtional} $K(t,f;X_0,X_1)\equiv K(t,f)$ is defined for $f\in X_0+X_1$ and $t>0$ by
 \begin{align*}
K(t,f)=\inf \Big \{\|f_0\|_{X_0}+t\|f_1\|_{X_1}:\ f=f_0+f_1,\ f_i\in X_i , \; i=0,1
\Big \}.
\end{align*}
It is well-known that for each $f\in X_0+X_1$, the function $t\rightsquigarrow K(t,f)$ is non-decreasing, while $t\rightsquigarrow t^{-1}K(t,f)$, $t>0$, is non-increasing. Other important property of the $K$-functional is the fact that 
\begin{equation}\label{eKK}
K(t,f;X_0,X_1)=tK(t^{-1},f;X_1,X_0)\quad \mbox{for all}\ t>0,
\end{equation}
(see \cite[Chap. 5, Proposition 1.2]{Bennett-Sharpley}).

We now recall the definition and some properties of the real interpolation method  $\overline{X}_{\theta,\b,E}$, the limiting $\mathcal{R}$ and $\mathcal{L}$ constructions 
and the interpolation functors in \eqref{extreme1}.
We refer to the papers \cite{FMS-1, FMS-RL1, FMS-RL2} for the major part of the results of this subsection, and for more information about these constructions.

\begin{defn}\label{defrealmethod}
Let $E$ be an r.i.  space,  $\b\in SV$ and $0\leq\theta \leq 1$. The real interpolation space $\overline{X}_{\theta,\b,E}\equiv(X_0,X_1)_{\theta, \b, E}$ consists of all $f$ in $X_{0} + X_{1}$ for which
$$ \|f\|_{\theta,\b,E} := \big \| t^{-\theta} {\b}(t) K(t,f) \big \|_{\widetilde{E}} < \infty.$$
\end{defn}

The space  $\overline{X}_{\theta, \b,E}$ is a (quasi-) Banach space, and it is an intermediate space for the couple $\overline{X}$, that is,
$$X_0\cap X_1\hookrightarrow \overline{X}_{\theta, \b,E}\hookrightarrow X_0+X_1,$$
provided that 
$0< \theta< 1$, or
 $\theta=0$ and $\|\b\|_{\widetilde{E}(1,\infty)}<\infty$ or
 $\theta=1$ and $\|\b\|_{\widetilde{E}(0,1)}<\infty$.
If none of these conditions holds, the space is trivial, that is $\overline{X}_{\theta, \b,E}=\{0\}$.

When $0<\theta<1$, $b\equiv 1$ and $E=L_q$, the space $\overline{X}_{\theta,b,E}$ coincides with the classical real interpolation space $\overline{X}_{\theta,q}$. When $\theta =0$ or $\theta =1$ this scale contains extrapolation spaces in the sense of Milman~\cite{Mi}, Gómez and Milman \cite{Go-Mi} and Astashkin, Lykov and Milman \cite{ALM}.  The spaces $\overline{X}_{\theta,\b,L_q}$ have been studied in detail by Gogatishvilli, Opic and Trebels in~\cite{GOT}. See also \cite{AEEK,clku,Do,EOP,gustav,per,tesisalba}, among other references.

\begin{rem}
\textup{Here we collect some elementary estimates. Using Lemma \ref{lem1} (i) and the monotonicity of $t\rightsquigarrow t^{-1}K(t,\cdot)$,  it is easy to check that for all $t>0$ and $f\in X_0+X_1$ 
\begin{equation}\label{e1}
t^{-\theta}\b(t)K(t,f)\lesssim\big \| s^{-\theta} {\b}(s) K(s,f) \big \|_{\widetilde{E}(0,t)},\qquad 0\leq\theta\leq1,\end{equation}
and
\begin{equation}\label{e3}
t^{-\theta}\b(t)K(t,f)\lesssim\big \| s^{-\theta} {\b}(s) K(s,f) \big \|_{\widetilde{E}(t,\infty)},
\qquad 0<\theta\leq 1.
\end{equation}}
\end{rem}

\begin{defn}
The (quasi-) Banach space 
$\overline{X}_{\theta,\b,E,a,F}^{\mathcal R}\equiv(X_0,X_1)_{\theta,\b,E,a,F}^{\mathcal R}$ consists of all $f\in X_0+X_1$ for which
$$\| f  \|_{\mathcal{R};\theta,\b,E,a,F} := \Big \|  \b(t) \|   s^{-\theta} a(s) K(s,f) \|_{\widetilde{F}(t,\infty)}      \Big   \|_{\widetilde{E}}<\infty.$$
\end{defn}

The space $\mathcal{R}$ is intermediate space for the couple $\overline{X}$, that is,
$$ X_0\cap X_1 \hookrightarrow  \overline{X}_{\theta,\b,E,a,F}^{\mathcal R} \hookrightarrow X_0+X_1$$
provided that any of the following conditions holds:
\begin{enumerate}
\item[1.] $0 < \theta < 1$ and  $\|\b\|_{\widetilde{E}(0,1)} <\infty$ or 
\item[2.] $\theta =0$,  $\|\b\|_{\widetilde{E}(0,1)} \!<\infty$ and   $\big \| \b(t)\|a\|_{\widetilde{F}(t,\infty)}  \big \|_{\widetilde{E}(1,\infty)} \!< \infty$
\item[3.] $\theta =1$, $ \|\b\|_{\widetilde{E}(0,1)} \!<\infty$,    $\big \| \b(t) \|a\|_{\widetilde{F}(t,1)} \big \|_{\widetilde{E}(0,1)} \!< \infty$ and  $\|\a\b\|_{\widetilde{E}(0,1)}<\infty$.
\end{enumerate}
Otherwise, $\overline{X}^{\mathcal R}_{\theta, \b,E,a,F}=\{0\}$.

\begin{defn}\label{defLR}
 Let $E$, $F$ be two r.i. spaces, $a, \b\in SV$ and $0\leq \theta\leq 1$. The space $\overline{X}_{\theta,\b,E,a,F}^{\mathcal L}\equiv(X_0,X_1)_{ \theta,\b,E,a,F}^{\mathcal L}$ consists of all $f \in X_{0} + X_{1}$ for which 
$$ \|f \|_{\mathcal L;\theta,\b,E,a,F}:=
\Big\|\b(t) \|s^{-\theta} a(s) K(s,f) \|_{\widetilde{F}(0,t)}\Big\|_{\widetilde{E}} < \infty. $$
\end{defn}
This is  a (quasi-) Banach space. Moreover it is intermediate for the couple $\overline{X}$, $$ X_0\cap X_1 \hookrightarrow  \overline{X}_{\theta,\b,E,a,F}^{\mathcal L} \hookrightarrow X_0+X_1,$$ provided that 
\begin{enumerate}
\item[1.] $0 < \theta < 1$  and  $\|\b\|_{\widetilde{E}(1,\infty)} <\infty$,  or  
\item[2.] $\theta =0$, $\|\b\|_{\widetilde{E}(1,\infty)} \!<\infty$, $ \big \| \b(t)\|a\|_{\widetilde{F}(1,t)}  \big \|_{\widetilde{E}(1,\infty)} \!< \infty$ and  $\|\a\b\|_{\widetilde{E}(1,\infty)}<\infty$ or 
\item[3.] $\theta =1$, $\|\b\|_{\widetilde{E}(1,\infty)} \!<\infty$ and $\big \| \b(t) \|a\|_{\widetilde{F}(0,t)} \big \|_{\widetilde{E}(0,1)} \!< \infty$.
\end{enumerate}
If none of these conditions holds, then $\overline{X}^{\mathcal L}_{\theta, \b,E,a,F}$ is  the trivial space.

The above definitions  generalize a previous family of spaces introduced by  Evans and Opic  \cite{EO} in which  $\a, \b$ are broken logarithms while $E, F$ remain within the classes $L_q$. Earlier versions of these spaces appeared in a paper  by Doktorskii \cite{Do}, with $a, \b$ powers of logarithms and  $\overline{X}$ an ordered couple.  These spaces also appear in  the work of  Gogatishvili, Opic and Trebels \cite{GOT} and Ahmed \textit{et al.} \cite{AEEK}. In all these cases the spaces $E$ and $F$ remain within the $L_q$ classes.

The spaces $\overline{X}_{\theta,\b,E}$, $\overline{X}^{\mathcal{R}}_{\theta,\b,E,a,F}$ and  $\overline{X}^{\mathcal{L}}_{\theta,\b,E,a,F}$ satisfy the following symmetry property.

\begin{lem}\label{symLR}
Let $E$, $F$ be r.i. spaces, $\a$, $\b\in SV$ and $0\leq\theta\leq1$. Then
$$(X_0,X_1)_{\theta,\b,E}=(X_1,X_0)_{1-\theta,\overline{\b},E}\mand(X_0,X_1)^{\mathcal{L}}_{\theta,\b,E,a,F}=(X_1,X_0)^{\mathcal{R}}_{1-\theta,\overline{\b},E,\overline{a},F}.$$
\end{lem}

Next lemma collects two estimates that will be used in the rest of the paper.

\begin{lem}\label{lemLRK}
Let $E$, $F$ be  r.i. spaces, $a, \b\in SV$ and $0\leq \theta\leq 1$. Then, for all $f\in X_0+X_1$ and $u>0$
\begin{equation}\label{e5}
u^{-\theta}\a(u)\|\b\|_{\widetilde{E}(0,u)}K(u,f)\lesssim \Big\| \b(t) \|s^{- \theta} \a(s) K(s,f)\|_{ \widetilde{F}(t,u)}\Big\|_{ \widetilde{E}(0,u)}
\end{equation}
and
\begin{equation}\label{e7}
u^{-\theta}\a(u)\|\b\|_{\widetilde{E}(u,\infty)}K(u,f)\lesssim \Big\| \b(t) \|s^{- \theta} \a(s) K(s,f)\|_{ \widetilde{F}(u,t)}\Big\|_{ \widetilde{E}(u,\infty)}.
\end{equation}
\end{lem}

\begin{proof}
We refer to \cite[Lemma 2.12]{FMS-RL1} for the proof \eqref{e5}. Next we prove \eqref{e7}. Let $f\in X_0+X_1$ and $u>0$. Using the monotonicity of the $K$-functional, Lemma \ref{lem0} (iii) and Lemma \ref{lem1} (ii) and (iv), we arrived at
\begin{align*}
\Big\| \b(t) \|s^{- \theta} \a(s) K(s,f)\|_{ \widetilde{F}(u,t)}\Big\|_{ \widetilde{E}(u,\infty)}&\geq
\Big\| \b(t) \|s^{- \theta} \a(s) K(s,f)\|_{ \widetilde{F}(u,t)}\Big\|_{ \widetilde{E}(2u,\infty)}\\
&\geq\Big\| \b(t) \|s^{- \theta} \a(s) K(s,f)\|_{ \widetilde{F}(u,2u)}\Big\|_{ \widetilde{E}(2u,\infty)}\\
&\sim u^{-\theta}\a(u)\|\b\|_{\widetilde{E}(2u,\infty)}K(u,f)\\
&\sim u^{-\theta}\a(u)\|\b\|_{\widetilde{E}(u,\infty)}K(u,f).
\end{align*}
\end{proof}

Consequently, when $0\leq\theta\leq1$ 
\begin{equation}\label{eK2}
K(u,f)\lesssim \frac{u^\theta}{\a(u)\|\b\|_{\widetilde{E}(0,u)}}\|f\|_{\mathcal R;\theta,\b,E,a,F}
\end{equation}
for all $u>0$ and $f\in\overline{X}^{\mathcal R}_{\theta,\b,E,\a,F}$, and
\begin{equation}\label{eK3}
K(u,f)\lesssim \frac{u^\theta}{\a(u)\|\b\|_{\widetilde{E}(u,\infty)}}\|f\|_{\mathcal L;\theta,\b,E,a,F}
\end{equation}
for all $u>0$ and $f\in\overline{X}^{\mathcal L}_{\theta,\b,E,\a,F}$.

The extreme classes in \eqref{extreme1} are defined in \cite{FMS-RL1,FMS-RL2} in the following way.

\begin{defn}\label{defLRR}
Let $E$, $F$, $G$ be  r.i. spaces, $\a, \b, c \in SV$ and $0< \theta<1$. The space
$\overline{X}_{\theta,c,E,\b,F,\a,G}^{\mathcal R,\mathcal R}\equiv(X_0,X_1)_{ \theta,c,E,\b,F,\a,G}^{\mathcal R,\mathcal R}$ is the set of all  $f\in X_0+X_1$ such that
\begin{equation}\label{dRR}
 \|f \|_{\mathcal R,\mathcal R;\theta,c,E,\b,F,\a,G} :=
\bigg\|c(u)\Big\|\b(t) \|s^{-\theta} \a(s) K(s,f) \|_{\widetilde{G}(t,\infty)}\Big\|_{\widetilde{F}(u,\infty)}\bigg\|_{\widetilde{E}} < \infty.
\end{equation}
The space
$\overline{X}_{\theta,c,E,\b,F,\a,G}^{\mathcal L,\mathcal L}\equiv(X_0,X_1)_{ \theta,c,E,\b,F,\a,G}^{\mathcal L,\mathcal L}$ is the set of all $f\in X_0+X_1$ for which 
\begin{equation}\label{dLL}
 \|f \|_{\mathcal L,\mathcal L;\theta,c,E,\b,F,\a,G} :=
\bigg\|c(u)\Big\|\b(t) \|s^{-\theta} \a(s) K(s,f) \|_{\widetilde{G}(0,t)}\Big\|_{\widetilde{F}(0,u)}\bigg\|_{\widetilde{E}}<\infty.
\end{equation}

The space $\overline{X}_{\theta,c,E,\b,F,\a,G}^{\mathcal R,\mathcal L}\equiv(X_0,X_1)_{ \theta,c,E,\b,F,\a,G}^{\mathcal R,\mathcal L}$ is the set of all $f\in X_0+X_1$ for which
\begin{equation}\label{dRL}
 \|f \|_{\mathcal R,\mathcal L;\theta,c,E,\b,F,a,G} :=
\bigg\|c(u)\Big\|\b(t) \|s^{-\theta} a(s) K(s,f) \|_{\widetilde{G}(t,u)}\Big\|_{\widetilde{F}(0,u)}\bigg\|_{\widetilde{E}}<\infty.
\end{equation}

The space $\overline{X}_{\theta,c,E,\b,F,a,G}^{\mathcal L,\mathcal R}\equiv(X_0,X_1)_{ \theta,c,E,\b,F,a,G}^{\mathcal L,\mathcal R}$ is the set of all $f\in X_0+X_1$ such that
\begin{equation}\label{dLR}
 \|f \|_{\mathcal L,\mathcal R;\theta,c,E,\b,F,\a,G} :=
\bigg\|c(u)\Big\|\b(t) \|s^{-\theta} a(s) K(s,f) \|_{\widetilde{G}(u,t)}\Big\|_{\widetilde{F}(u,\infty)}\bigg\|_{\widetilde{E}}<\infty.
\end{equation}
\end{defn}

These spaces enjoy the following symmetry property.
\begin{lem}\label{symLL}
Let $E$, $F$, $G$ be r.i. spaces, $a$, $\b$, $c\in SV$ and $0\leq\theta\leq1$. Then
$$(X_0,X_1)^{\mathcal{L,L}}_{\theta,c,E,\b,F,a,G}=(X_1,X_0)^{\mathcal{R,R}}_{1-\theta,\overline{c},E,\overline{b},F,\overline{a},G}$$
and
$$(X_0,X_1)^{\mathcal{L,R}}_{\theta,c,E,\b,F,a,G}=(X_1,X_0)^{\mathcal{R,L}}_{1-\theta,\overline{c},E,\overline{b},F,\overline{a},G}.$$
\end{lem}

\section{Generalized Holmstedt type formulae}\label{sectionH}
For  parameters $0<\theta_0<\theta_1< 1$, $a_0$, $a_1$, $\b_0$, $\b_1\in SV$ and $E_0$, $E_1$, $F_0$, $F_1$ r.i spaces, the couples
$$(\overline{X}^{\mathcal R}_{\theta_0,\b_0,E_0,\a_0,F_0}, \overline{X}^{\mathcal R}_{\theta_1, \b_1,E_1,\a_1,F_1}),\quad
(\overline{X}^{\mathcal L}_{\theta_0, \b_0,E_0,\a_0,F_0}, \overline{X}^{\mathcal L}_{\theta_1,\b_1,E_1,\a_1,F_1})$$
and
$$ (\overline{X}^{\mathcal R}_{\theta_0,\b_0,E_0,\a_0,F_0}, \overline{X}^{\mathcal L}_{\theta_1, \b_1,E_1,\a_1,F_1}),\quad
(\overline{X}^{\mathcal L}_{\theta_0, \b_0,E_0,\a_0,F_0}, \overline{X}^{\mathcal R}_{\theta_1,\b_1,E_1,\a_1,F_1})$$
are compatible (quasi-) Banach couples. We wish to identify the interpolation spaces generated by the previous couples as interpolated spaces of the original one, $\overline{X}$. Our first stage in that process is to relate the $K$-functional of the underlying couples through several generalized Holmstedt type formulae. 

\subsection{The $K$-functional of the couple $(\overline{X}^{\mathcal R}_{\theta_0,\b_0,E_0,\a_0,F_0},\overline{X}^{\mathcal R}_{\theta_1,\b_1,E_1,\a_1,F_1})$, $0<\theta_0<\theta_1<1$.}
\begin{thm} \label{5.12}
Let $0<\theta_0<\theta_1<1$. Let $E_0$, $E_1$, $F_0$, $F_1$ r.i. spaces and $\a_0$, $\a_1$, 
$\b_0$,  $\b_1\in SV$ such that $\|\b_0\|_{\widetilde{E}_0(0,1)}<\infty$ and $\|\b_1\|_{\widetilde{E}_1(0,1)}<\infty$.  Then, for every $f \in \overline{X}^{\mathcal R}_{\theta_0, \b_0,E_0,\a_0,F_0}+\overline{X}^{\mathcal R}_{\theta_1, \b_1,E_1,\a_1,F_1} $ and $u>0$
 \begin{align}
K\big( \rho(u), f; \overline{X}^{\mathcal R}_{\theta_0, \b_0,E_0,\a_0,F_0},
&\overline{X}^{\mathcal R}_{\theta_1, \b_1,E_1,\a_1,F_1}\big)\label{eqHoRR} \\
 &\sim\Big\| \b_0(t) \|s^{- \theta_0} \a_0(s) K(s,f)\|_{ \widetilde{F}_0(t,u)}\Big\|_{ \widetilde{E}_0(0,u)}\nonumber \\
& + \rho(u)\| \b_1\|_{ \widetilde{E}_1(0,u)} \cdot\|t^{- \theta_1} \a_1(t) K(t,f)\|_{ \widetilde{F}_1(u,\infty)}\nonumber\\
&+\rho(u)\Big\| \b_1(t) \|s^{- \theta_1} \a_1(s) K(s,f)\|_{ \widetilde{F}_1(t,\infty)}\Big\|_{ \widetilde{E}_1(u,\infty)}\nonumber
\end{align}
where
\begin{equation}\label{rhoRR}
\rho(u) =u^{\theta_1- \theta_0} \frac{\a_0(u)\|\b_0\|_{\widetilde{E}_0(0,u)}}{\a_1(u)\|\b_1\|_{\widetilde{E}_1(0,u)} },\quad u >0.
\end{equation}
\end{thm}

\begin{proof}
Given $f\in X_0+X_1$ and $u>0$, we consider the following (quasi-) norms 
\begin{align}
(P_j f)(u) &= \Big\| \b_j(t) \|s^{- \theta_j} \a_j(s) K(s,f)\|_{ \widetilde{F}_j(t,u)}\Big\|_{ \widetilde{E}_j(0,u)},\nonumber \\
(R_j f)(u) &= \| \b_j\|_{ \widetilde{E}_j(0,u)} \cdot\|t^{- \theta_j} \a_j(t) K(t,f)\|_{ \widetilde{F}_j(u,\infty)},\label{epqr}\\
(Q_j f)(u) &= \Big\| \b_j(t)\| s^{- \theta_j} \a_j(s) K(s,f)\|_{ \widetilde{F}_j(t,\infty)}\Big\|_{\widetilde{E}_j(u,\infty)},\nonumber
\end{align}
and we denote $Y_j= \overline{X}^{\mathcal R}_{\theta_j, \b_j,E_j,\a_j,F_j}$, $j=0, 1$. With this notation what we pursue to show is  that
\begin{equation}\label{5173}
K(\rho(u), f;Y_0,Y_1) \sim (P_{0}f)(u) + \rho(u) [(R_1f)(u)+(Q_{1}f)(u)]
\end{equation}
where $\rho$ is defined by  \eqref{rhoRR}.

We fix $f\in X_0+X_1$, $u>0$ and we assume that $(P_{0}f)(u)$, $(R_1f)(u)$ and $(Q_{1}f)(u)$ are finite, otherwise the upper estimate of  \eqref{5173} holds trivially. 
As usual (see for example~\cite{Bennett-Sharpley} or \cite{EOP}) we choose a decomposition $f = g + h$
such that $ \|g\|_{X_0} + u\|h\|_{X_1} \leq 2 K(u,f)$. Then, for all $s>0$,
\begin{equation}\label{4552}
K(s,g) \leq 2 K(u,f)\ \ \text{ and }  \ \ K(s,h) \leq 2s \frac{K(u,f)}{u}.
\end{equation}
Therefore, in order to obtain the upper estimate of (\ref{5173}) it is enough to prove that
$$
\|g\|_{Y_0}+\rho(u)\|h\|_{Y_1}\lesssim (P_0f)(u)+\rho(u) [(R_1f)(u)+(Q_{1}f)(u)].
$$
We start establishing that $\|g\|_{Y_0}\lesssim (P_0 f)(u)$. 
By the triangle inequality and the (quasi-) subadditivity of the $K$-functional, we have the inequalities
\begin{align}\label{456}
\|g\|_{Y_0}&\leq(P_0g)(u)+(R_0g)(u)+(Q_0g)(u)\\&\lesssim (P_0f)(u)+(P_0h)(u)+(R_0g)(u)+(Q_0g)(u).\nonumber
\end{align}
So, it suffices to estimate $(P_0 h)(u)$, $(R_0g)(u)$ and $(Q_0 g)(u) $ from above. Using (\ref{4552}) and the monotonicity of $s\rightsquigarrow s^{-1}K(s,f)$, $s>0$,  we  obtain that
\begin{align}
(P_0 h)(u) &= \Big\| \b_0(t) \|s^{- \theta_0} \a_0(s) K(s,h)\|_{ \widetilde{F}_0(t,u)}\Big\|_{ \widetilde{E}_0(0,u)}\label{457}\\
&\lesssim \Big\| \b_0(t) \|s^{1- \theta_0} \a_0(s)\frac{K(u,f)}{u}\|_{ \widetilde{F}_0(t,u)}\Big\|_{ \widetilde{E}_0(0,u)}\nonumber\\
&\leq \Big\| \b_0(t) \|s^{- \theta_0} \a_0(s)K(s,f)\|_{ \widetilde{F}_0(t,u)}\Big\|_{ \widetilde{E}_0(0,u)}=(P_0 f)(u).\nonumber
\end{align}
On the other hand, using  (\ref{4552}), Lemma \ref{lem1} (i) and \eqref{e5}, we have
\begin{align}
(R_0 g)(u) &= \| \b_0\|_{ \widetilde{E}_0(0,u)} \|t^{- \theta_0} \a_0(t) K(t,g)\|_{ \widetilde{F}_0(u,\infty)}\label{458}\\
&\lesssim K(u,f) \| \b_0\|_{ \widetilde{E}_0(0,u)} \|t^{- \theta_0} \a_0(t) \|_{ \widetilde{F}_0(u,\infty)}\nonumber\\
&\sim u^{-\theta_0}\a_0(u)\|\b_0\|_{\widetilde{E}_0(0,u)}K(u,f)\lesssim (P_0 f)(u).\nonumber
\end{align}
Similarly, using also Lemma \ref{lem1} (iv) it follows that
\begin{align}
(Q_0 g)(u) &= \Big\| \b_0(t) \|s^{- \theta_0} \a_0(s) K(s,g)\|_{ \widetilde{F}_0(t,\infty)}\Big\|_{ \widetilde{E}_0(u,\infty)}\label{459}\\
&\lesssim K(u,f)\Big\| \b_0(t) \|s^{- \theta_0} \a_0(s)\|_{ \widetilde{F}_0(t,\infty)}\Big\|_{ \widetilde{E}_0(u,\infty)}\nonumber\\
&\sim u^{-\theta_0}\b_0(u)\a_0(u)K(u,f)\lesssim (P_0 f)(u).\nonumber
\end{align}
Hence, inequalities (\ref{456})-(\ref{459}) give  that $\|g\|_{Y_0}\lesssim (P_0 f)(u)<\infty$, which in particular shows that $g\in Y_0$. In a similar way it can be proved that
\begin{equation*}\label{461}
\|h\|_{Y_1}\lesssim(R_1f)(u)+(Q_{1}f)(u)<\infty,
\end{equation*}
so $h\in Y_1$. Indeed, by  the triangle inequality and the (quasi-) subadditivity of the $K$-functional, we have
\begin{align}
\|h\|_{Y_1}&\leq (P_1h)(u)+(R_1h)(u)+(Q_1h)(u)\label{463}\\
&\lesssim (P_1h)(u)+(R_1f)(u)+(R_1g)(u)+(Q_1f)(u)+(Q_1g)(u).\nonumber
\end{align}
Now we estimate $(P_1h)(u)$. Using (\ref{4552}), Lemma \ref{lem1} (i) and \eqref{e3},  we obtain
\begin{align}
(P_1 h)(u)&=\Big\| \b_1(t) \|s^{- \theta_1} \a_1(s) K(s,h)\|_{ \widetilde{F}_1(t,u)}\Big\|_{ \widetilde{E}_1(0,u)} \label{464}\\
&\lesssim \frac{K(u,f)}{u}\Big\| \b_1(t) \|s^{1- \theta_1} \a_1(s)\|_{ \widetilde{F}_1(0,u)}\Big\|_{ \widetilde{E}_1(0,u)} \nonumber\\
&\sim u^{- \theta_1} \a_1(u)\| \b_1\|_{\widetilde{E}_1(0,u)}  K(u,f)\lesssim (R_1 f)(u).\nonumber
 \end{align}
Similarly,  $(R_1 g)(u)$
and
$(Q_1g)(u)$ can be bounded from above. In fact, Lemma \eqref{4552}, Lemma \ref{lem1} (i)  and \eqref{e3} yield
\begin{align}
(R_1 g)(u)&=\|\b_1\|_{\widetilde{E}_1(0,u)} \big\|t^{- \theta_1} \a_1(t) K(t,g)\big\|_{ \widetilde{F}_1(u,\infty)}\label{465} \\
&\lesssim K(u,f)\|\b_1\|_{\widetilde{E}_1(0,u)} \big\|t^{- \theta_1} \a_1(t)\big\|_{ \widetilde{F}_1(u,\infty)} \nonumber\\
&\sim u^{-\theta_1}\a_1(u)\|\b_1\|_{\widetilde{E}_1(0,u)} K(u,f)\lesssim (R_1 f)(u),\nonumber
\end{align}
and also, using Lemma \ref{lem1} (iv),
 \begin{align}
(Q_1 g)(u)&=\Big\| \b_1(t) \|s^{- \theta_1} \a_1(s) K(s,g)\|_{ \widetilde{F}_1(t,\infty)}\Big\|_{ \widetilde{E}_1(u,\infty)}\label{466} \\
&\lesssim K(u,f)\Big\| \b_1(t) \|s^{- \theta_1} \a_1(s)\|_{ \widetilde{F}_1(t,\infty)}\Big\|_{ \widetilde{E}_1(u,\infty)}\nonumber \\
&\lesssim u^{- \theta_1} \a_1(u)\b_1(u) K(u,f)\lesssim (R_1 f)(u).\nonumber
 \end{align}
Summing up
$$\|g\|_{Y_0}+\rho(u)\|h\|_{Y_1}\lesssim (P_0f)(u)+\rho(u)[(R_1 f)(u)+(Q_{1}f)(u)],$$
for any positive function $\rho$. This concludes the proof of the upper estimate (\ref{5173}).

Next we proceed with the lower estimate of \eqref{5173}, that is 
$$(P_0f)(u)+\rho(u)[(R_1 f)(u)+(Q_1f)(u)]\lesssim K(\rho(u), f;Y_0,Y_1).$$
Suposse that $f\in Y_0+Y_1$ and fix $u>0$. Let $f=f_0+f_1$ be any decomposition of $f$ with $f_0\in  Y_0$ and $f_1\in Y_1$. The (quasi-) subadditivity of the $K$-functional and the definition of the norm in $Y_0$, $Y_1$ imply\begin{align*}
(P_{0}f)(u) & \lesssim  (P_{0}f_0)(u) + (P_{0}f_1)(u)  \leq  \|f_0\|_{Y_{0}} + (P_{0}f_1)(u), \\
(R_{1}f)(u) & \lesssim (R_{1}f_0)(u) + (R_{1}f_1)(u)  \leq
(R_{1}f_0)(u) + \|f_1\|_{Y_{1}},\\
(Q_{1}f)(u) & \lesssim (Q_{1}f_0)(u) + (Q_{1}f_1)(u)  \leq
(Q_{1}f_0)(u) + \|f_1\|_{Y_{1}}.
\end{align*}
Then
\begin{align*}
(P_0f)(u)+&\rho(u)[(R_1 f)(u)+(Q_1f)(u)]\\
&\lesssim  \|f_0\|_{Y_{0}} + \rho(u)\|f_1\|_{Y_1}+ (P_{0}f_1)(u) +\rho(u)[(R_{1}f_0)(u)+(Q_{1}f_0)(u)].
\end{align*}
Thus, to finish the proof, it is enough to verify that
\begin{align*}
(P_{0}f_1)(u) + \rho(u) [(R_1 f_0)(u)+(Q_{1}f_0)(u)]\lesssim \| f_0 \|_{Y_{0}} + \rho(u)\|f_1\|_{Y_{1}} .
\end{align*}
We begin  with $(P_{0}f_1)(u)$.  Using \eqref{eK2} with $f=f_1\in Y_1$ and Lemma \ref{lem1} (i) ($\theta_1-\theta_0>0$), we obtain 
\begin{align}
(P_0f_1)(u)
&\lesssim \|f_1\|_{Y_1}\bigg\| \b_0(t) \Big\|s^{\theta_1- \theta_0} \frac{\a_0(s)}{\a_1(s)\|\b_1\|_{\widetilde{E}_1(0,s)}}\Big\|_{ \widetilde{F}_0(0,u)}\bigg\|_{ \widetilde{E}_0(0,u)}\label{462}\\
&\sim u^{\theta_1- \theta_0} \frac{\a_0(u)\|\b_0\|_{ \widetilde{E}_0(0,u)}}{\a_1(u)\|\b_1\|_{\widetilde{E}_1(0,u)}} \|f_1\|_{Y_1}=\rho(u)\|f_1\|_{Y_1}.\nonumber
\end{align}
In a similar way, using \eqref{eK2} with $f=f_0\in Y_0$ and Lemma \ref{lem1} (i) ($\theta_0-\theta_1<0$), we have 
\begin{align}
(R_1f_0)(u)&\lesssim \|f_0\|_{Y_0}\| \b_1\|_{\widetilde{E}_1(0,u)}\Big\|t^{\theta_0- \theta_1} \frac{\a_1(t)}{\a_0(t)\|\b_0\|_{\widetilde{E}_0(0,t)}}\Big\|_{ \widetilde{F}_1(u,\infty)}\label{463bis}\\
&\sim u^{\theta_0-\theta_1}\frac{\a_1(u)\|\b_1\|_{\widetilde{E}_1(0,u)}}{\a_0(u)\|\b_0\|_{\widetilde{E}_0(0,u)}}\|f_0\|_{Y_0}=\frac1{\rho(u)}\|f_0\|_{Y_0}\nonumber
\end{align}
and
\begin{align}
(Q_1f_0)(u)
&\lesssim \|f_0\|_{Y_0}\bigg\| \b_1(t)\Big\|s^{\theta_0- \theta_1} \frac{\a_1(s)}{\a_0(s)\|\b_0\|_{\widetilde{E}_0(0,s)}}\Big\|_{ \widetilde{F}_1(t,\infty)}\bigg\|_{ \widetilde{E}_1(u,\infty)}\label{464bis}\\
&\sim u^{\theta_0-\theta_1}\frac{\a_1(u)\b_1(u)}{\a_0(u)\|\b_0\|_{\widetilde{E}_0(0,u)}}\|f_0\|_{Y_0}\lesssim\frac1{\rho(u)}\|f_0\|_{Y_0},\nonumber
\end{align}
where the last equivalence follows from Lemma \ref{lem1} iv).

Putting together the previous estimates  we establish that
$$(P_{0}f)(u) + \rho(u)[(R_1 f)(u)+ (Q_{1}f)(u)\lesssim \|f_0\|_{Y_0}+\rho(u)\|f_1\|_{Y_1}.$$
And taking infimum over all possible decomposition of $f=f_0+f_1$, with $f_0\in Y_0$ and $f_1\in Y_1$, we finish the proof.
\end{proof}

\subsection{The $K$-functional of the couple $(\overline{X}^{\mathcal L}_{\theta_0,\b_0,E_0,\a_0,F_0},\overline{X}^{\mathcal L}_{\theta_1,\b_1,E_1,\a_1,F_1})$, $0<\theta_0<\theta_1<1$.}\hspace*{1mm}

Next theorem follow from Theorem \ref{5.12} by means of a symmetry argument.

\begin{thm} \label{5.13}
Let $0<\theta_0<\theta_1<1$. Let $E_0$, $E_1$, $F_0$, $F_1$ r.i. spaces and $\a_0$, $\a_1$, $\b_0$, $\b_1\in SV$ such that $\|\b_0\|_{\widetilde{E}_0(1,\infty)}<\infty$ and $\|\b_1\|_{\widetilde{E}_1(1,\infty)}<\infty$.
Then, for every $f \in \overline{X}^{\mathcal L}_{\theta_0, \b_0,E_0,\a_0,F_0}+\overline{X}^{\mathcal L}_{\theta_1, \b_1,E_1,\a_1,F_1} $ and $u>0$
\vspace{2mm}
 \begin{align}
K\big( \rho(u), f; \overline{X}^{\mathcal L}_{\theta_0, \b_0,E_0,\a_0,F_0},
&\overline{X}^{\mathcal L}_{\theta_1, \b_1,E_1,\a_1,F_1}\big) \label{eHL}\\
 &\sim\Big\| \b_0(t) \|s^{- \theta_0} \a_0(s) K(s,f)\|_{ \widetilde{F}_0(0,t)}\Big\|_{ \widetilde{E}_0(0,u)} \nonumber\\
 &+\| \b_0\|_{\widetilde{E}_0(u,\infty)} \|t^{- \theta_0} \a_0(t) K(t,f)\|_{ \widetilde{F}_0(0,u)}\nonumber\\
& +\rho(u)\Big\| \b_1(t) \|s^{- \theta_1} \a_1(s) K(s,f)\|_{ \widetilde{F}_1(u,t)}\Big\|_{ \widetilde{E}_1(u,\infty)}\nonumber
\end{align}
where
\begin{equation*}\label{rhoLL}
\rho(u) =u^{\theta_1- \theta_0} \frac{\a_0(u)\|\b_0\|_{\widetilde{E}_0(u,\infty)}}{\a_1(u)\|\b_1\|_{\widetilde{E}_1(u,\infty)}},\quad  u>0.
\end{equation*}
\end{thm}

\begin{proof}
Lemma \ref{symLR} and \eqref{eKK} yield
\begin{align*}
K\big(\rho(u), f;&\overline{X}^{\mathcal L}_{\theta_0,\b_0,E_0,\a_0,F_0},\overline{X}^{\mathcal L}_{\theta_1, \b_1,E_1,\a_1,F_1}\big)\\
&=\rho(u)K\Big(\frac1{\rho(u)},f;(X_1,X_0)^{\mathcal R}_{1-\theta_1, \overline{\b}_1,E_1,\overline{\a}_1,F_1},(X_1,X_0)^{\mathcal R}_{1-\theta_0,\overline{\b}_0,E_0,\overline{\a}_0,F_0}\Big),
\end{align*}
where
$$\frac{1}{\rho(u)}=u^{\theta_0-\theta_1}\frac{\overline{\a}_1(\frac1u)}{\overline{\a}_0(\frac1u)}\frac{\|\overline{\b}_1\|_{\widetilde{E}_1(0,1/u)}}{\|\overline{\b}_0\|_{\widetilde{E}_0(0,1/u)}},\quad u>0.$$
Now, applying Theorem \ref{5.12} we obtain the estimate
\begin{align*}
 K\big(\rho(u), f;\overline{X}^{\mathcal L}_{\theta_0,\b_0,E_0,a_0,F_0},&\overline{X}^{\mathcal L}_{\theta_1, \b_1,E_1,a_1,F_1}\big)\\
 &\sim\rho(u)\Big\| \overline{\b}_1(t) \|s^{\theta_1-1} \overline{a}_1(s) K(s,f;X_1,X_0)\|_{ \widetilde{F}_1(t,1/u)}\Big\|_{ \widetilde{E}_1(0,1/u)} \\
& +\|\overline{\b}_0\|_{ \widetilde{E}_0(0,1/u)}\|t^{ \theta_0-1}\overline{a}_0(t) K(t,f;X_1,X_0)\|_{ \widetilde{F}_0(1/u,\infty)}\\
&+\Big\| \overline{\b}_0(t) \|s^{\theta_0-1} \overline{a}_0(s) K(s,f;X_1,X_0)\|_{ \widetilde{F}_0(t,\infty)}\Big\|_{ \widetilde{E}_0(1/u,\infty)}.
\end{align*}
The equivalence \eqref{eHL} follows using again \eqref{fbar} and \eqref{eKK}.
\end{proof}

\subsection{The $K$-functional of the couple $(\overline{X}^{\mathcal R}_{\theta_0,\b_0,E_0,\a_0,F_0},\overline{X}^{\mathcal L}_{\theta_1,\b_1,E_1,\a_1,F_1})$, $0<\theta_0<\theta_1<1$.}

\begin{thm} \label{5.17}
Let $0<\theta_0<\theta_1<1$. Let $E_0$, $E_1$, $F_0$, $F_1$ r.i. spaces and $\a_0$, $\a_1$, $\b_0$, $\b_1\in SV$ such that $\|\b_0\|_{\widetilde{E}_0(0,1)}<\infty$ and $\|\b_1\|_{\widetilde{E}_1(1,\infty)}<\infty$. 
Then, for every $f \in \overline{X}^{\mathcal R}_{\theta_0, \b_0,E_0,\a_0,F_0}+\overline{X}^{\mathcal L}_{\theta_1, \b_1,E_1,\a_1,F_1} $ and $u>0$
\begin{align*}
K\big( \rho(u), f; \overline{X}^{\mathcal R}_{\theta_0, \b_0,E_0,\a_0,F_0},
&\overline{X}^{\mathcal L}_{\theta_1, \b_1,E_1,\a_1,F_1}\big)\label{KRL} \\
 &\sim\Big\| \b_0(t) \|s^{- \theta_0} \a_0(s) K(s,f)\|_{ \widetilde{F}_0(t,u)}\Big\|_{ \widetilde{E}_0(0,u)} \nonumber\\
 &+\rho(u) \Big\| \b_1(t)\| s^{- \theta_1} \a_1(s) K(s,f)\|_{ \widetilde{F}_1(u,t)}\Big\|_{\widetilde{E}_1(u,\infty)},\nonumber
\end{align*}
where
\begin{equation}\label{rhoRL}
\rho(u) =u^{\theta_1- \theta_0} \frac{\a_0(u)\|\b_0\|_{\widetilde{E}_0(0,u)}}{\a_1(u)\|\b_1\|_{\widetilde{E}_1(u,\infty)} },\quad u >0.
\end{equation}
\end{thm}

\begin{proof} 
Given $f\in X_0+X_1$ and $u>0$, we consider the (quasi-) norms $(P_0 f)(u)$, $(R_0 f)(u)$ and $(Q_0 f)(u)$ defined as in \eqref{epqr} and we redefine
\begin{align*}
(P_1 f)(u) &= \Big\| \b_1(t) \|s^{- \theta_1} \a_1(s) K(s,f)\|_{ \widetilde{F}_1(0,t)}\Big\|_{ \widetilde{E}_1(0,u)} ,\\
(R_1 f)(u) &= \|\b_1\|_{ \widetilde{E}_1(u,\infty)} \cdot\|t^{- \theta_1} \a_1(t) K(t,f)\|_{ \widetilde{F}_1(0,u)},\\
(Q_1 f)(u) &= \Big\| \b_1(t)\| s^{- \theta_1} \a_1(s) K(s,f)\|_{ \widetilde{F}_1(u,t)}\Big\|_{\widetilde{E}_1(u,\infty)}.
\end{align*}
We denote, as usual, $Y_0= \overline{X}^{\mathcal R}_{\theta_0, \b_0,E_0,\a_0,F_0}$ and  $Y_1= \overline{X}^{\mathcal L}_{\theta_1, \b_1,E_1,\a_1,F_1}$.
What we want to show is that
\begin{equation}\label{5.17.3.3}
K(\rho(u), f;Y_0,Y_1) \sim (P_{0}f)(u) + \rho(u)(Q_{1}f)(u)
\end{equation}
where $\rho$ is defined by \eqref{rhoRL}. In order to do that we follow the same steps as in the proof of Theorem \ref{5.12}.

We fix $f \in X_0 + X_1$ and $u>0$, and we  choose a decomposition of $f = g + h$
such that $ \|g\|_{X_0} + u\|h\|_{X_1} \leq 2 K(u,f)$ so that \eqref{4552} is satisfied. We start by showing that 
$$\|g\|_{X_0} + \rho(u)\|h\|_{X_1}\lesssim (P_{0}f)(u) + \rho(u)(Q_{1}f)(u)$$
to establish the upper bound for the $K$-functional in \eqref{5.17.3.3}. The inequality $\|g\|_{Y_0}\lesssim (P_0f)(u)$ can be proved exactly as we did through \eqref{456} to \eqref{459}. 
Now we proceed with $\|h\|_{Y_1}$. Using the triangular inequality and the (quasi-) subadditivity of the $K$-functional it follows
\begin{align*}
\|h\|_{Y_1}&\leq (P_1h)(u)+(R_1h)(u)+(Q_1h)(u)\\
&\lesssim (P_1h)(u)+(R_1h)(u)+(Q_1f)(u)+(Q_1g)(u).
\end{align*}
Hence, it is enough to estimate $(P_{1}h)(u)$, $(R_1h)(u)$ and $(Q_{1}g)(u)$ from above. Using \eqref{4552}, Lemma \ref{lem1} (i), (iv) and \eqref{e7}, we deduce that
\begin{align*}
(P_1 h)(u)&=\Big\| \b_1(t) \|s^{- \theta_1} \a_1(s) K(s,h)\|_{ \widetilde{F}_1(0,t)}\Big\|_{ \widetilde{E}_1(0,u)} \\
&\lesssim \frac{K(u,f)}{u}\Big\| \b_1(t) \|s^{1- \theta_1} \a_1(s)\|_{ \widetilde{F}_1(0,t)}\Big\|_{ \widetilde{E}_1(0,u)} \\
&\sim u^{- \theta_1} \a_1(u)\b_1(u) K(u,f)\lesssim (Q_1 f)(u)
 \end{align*}
and
 \begin{align*}
(R_{1}h)(u) &=\|\b_1\|_{ \widetilde{E}_1(u,\infty)} \|t^{- \theta_1} \a_1(t) K(t,h)\|_{ \widetilde{F}_1(0,u)}\\
&\lesssim \frac{K(u,f)}{u}\|\b_1\|_{ \widetilde{E}_1(u,\infty)} \|t^{1- \theta_1} \a_1(t) \|_{ \widetilde{F}_1(0,u)}\\
&\sim u^{- \theta_1} \a_1(u)\|\b_1\|_{ \widetilde{E}_1(u,\infty)}  K(u,f)\lesssim (Q_1 f)(u).
\end{align*}
In the same vein
 \begin{align*}
(Q_{1}g)(u) &= \Big\| \b_1(t)\| s^{- \theta_1} \a_1(s) K(s,g)\|_{ \widetilde{F}_1(u,t)}\Big\|_{\widetilde{E}_1(u,\infty)}\\
&\lesssim K(u,f) \|\b_1\|_{ \widetilde{E}_1(u,\infty)}\| s^{- \theta_1} \a_1(s) \|_{ \widetilde{F}_1(u,\infty)}\\
&\sim u^{- \theta_1} \a_1(u)\|\b_1\|_{ \widetilde{E}_1(u,\infty)}  K(u,f)\lesssim (Q_1 f)(u).
 \end{align*}
Then,  $\|h\|_{Y_1}\lesssim (Q_{1}f)(u)$ and $h\in Y_1$. 
Summing up, we establish the  upper estimate of (\ref{5.17.3.3}), that is
$$K(\rho(u), f;Y_0,Y_1)\leq \|g\|_{Y_0}+\rho(u)\|h\|_{Y_1}\lesssim (P_0f)(u)+\rho(u)(Q_{1}f)(u).$$

Let us prove the lower estimate
\begin{equation}\label{860}
(P_0f)(u)+\rho(u)(Q_1f)(u)\lesssim K(\rho(u), f;Y_0,Y_1)
\end{equation}
for any $f\in Y_0+Y_1$ and $u>0$. 
Let $f=f_0+f_1$  be any decomposition of $f$ with $f_0\in Y_0$ and $f_1\in Y_1$. Using the (quasi-) subadditivity of the $K$-functional and the definition of the norm in $Y_0$, $Y_1$, we have
\begin{align*}
(P_{0}f)(u) & \lesssim(P_{0}f_0)(u) + (P_{0}f_1)(u)  \leq  \|f_0\|_{Y_{0}} + (P_{0}f_1)(u) \\
(Q_{1}f)(u) & \lesssim (Q_{1}f_0)(u) + (Q_{1}f_1)(u)  \leq
(Q_{1}f_0)(u) + \|f_1\|_{Y_{1}}.
\end{align*}
Thus, we have to study the boundedness of $(P_0 f_1)(u)$ and $\rho(u)(Q_1 f_0)(u)$  by $\| f_0 \|_{Y_{0}} + \rho(u)\|f_1\|_{Y_{1}}$. For the proof of the estimate $(P_0 f_1)(u)\lesssim \rho(u)\|f_1\|_{Y_{1}}$ one has to argue as in \eqref{462}.
Similarly, using \eqref{eK3} with $f=f_0$ and Lemma \ref{lem1} (i) ($\theta_0-\theta_1<0$), we obtain 
\begin{align*}
(Q_1f_0)(u)
&\lesssim \|f_0\|_{Y_0}\| \b_1\|_{\widetilde{E}_1(u,\infty)}\Big\| s^{\theta_0- \theta_1} \frac{\a_1(s)}{\a_0(s)\|\b_0\|_{\widetilde{E}_0(0,s)}}\Big\|_{ \widetilde{F}_1(u,\infty)}\\
&\sim u^{\theta_0- \theta_1} \frac{\a_1(u)\| \b_1\|_{\widetilde{E}_1(u,\infty)}}{\a_0(u)\|\b_0\|_{\widetilde{E}_0(0,u)}}\|f_0\|_{Y_0}=\frac{1}{\rho(u)}\|f_0\|_{Y_0}.\label{ec38}
\end{align*}
Putting together the previous estimates we establish that
$$(P_{0}f)(u) + \rho(u) (Q_{1}f)(u)\lesssim \|f_0\|_{Y_0}+\rho(u)\|f_1\|_{Y_1}.$$
Taking infimum over all possible decomposition of $f=f_0+f_1$, with $f_0\in Y_0$ and $f_1\in Y_1$, we deduce \eqref{860}. The proof is finished. 
\end{proof}

\subsection{The $K$-functional of the couple $(\overline{X}^{\mathcal L}_{\theta_0,\b_0,E_0,\a_0,F_0},\overline{X}^{\mathcal R}_{\theta_1,\b_1,E_1,\a_1,F_1})$, $0<\theta_0<\theta_1<1$.}\hspace*{1mm}

Finally we present the last Holmstedt type formula. Although it seems that the result is the symmetric counterpart of Theorem \ref{5.17}, this is not the case since the condition $\theta_0<\theta_1$ is crucial.

\begin{thm} \label{5.17.2}
Let $0<\theta_0<\theta_1<1$. Let $E_0$, $E_1$, $F_0$, $F_1$ r.i. spaces and $\a_0$, $\a_1$, $\b_0$, $\b_1\in SV$ such that $\|\b_0\|_{\widetilde{E}_0(1,\infty)}<\infty$ and $\|\b_1\|_{\widetilde{E}_1(0,1)}<\infty$. 
Then, for every $f \in \overline{X}^{\mathcal L}_{\theta_0, \b_0,E_0,\a_0,F_0}+\overline{X}^{\mathcal R}_{\theta_1, \b_1,E_1,\a_1,F_1} $ and $u>0$
 \begin{align*}
K\big( \rho(u), f; \overline{X}^{\mathcal L}_{\theta_0, \b_0,E_0,\a_0,F_0},
&\overline{X}^{\mathcal R}_{\theta_1, \b_1,E_1,\a_1,F_1}\big)\\
 &\sim\Big\| \b_0(t) \|s^{- \theta_0} \a_0(s) K(s,f)\|_{ \widetilde{F}_0(0,t)}\Big\|_{ \widetilde{E}_0(0,u)}\nonumber \\
 &+\|\b_0\|_{\widetilde{E}_0(u,\infty) }\|t^{- \theta_0} \a_0(t) K(t,f)\|_{ \widetilde{F}_0(0,u)}
\nonumber\\
&+\rho(u)  \|\b_1\|_{\widetilde{E}_1(0,u)} \|t^{- \theta_1} \a_1(t) K(t,f)\|_{ \widetilde{F}_1(u,\infty)}\nonumber\\
&+\rho(u) \Big\| \b_1(t)\| s^{- \theta_1} \a_1(s) K(s,f)\|_{ \widetilde{F}_1(t,\infty)}\Big\|_{\widetilde{E}_1(u,\infty)},\nonumber
\end{align*}
where
\begin{equation}\label{rhoLR}
\rho(u) =u^{\theta_1- \theta_0} \frac{\a_0(u)\|\b_0\|_{\widetilde{E}_0(u,\infty)}}{\a_1(u)\|\b_1\|_{\widetilde{E}_1(0,u)} },\quad u >0.\end{equation} 
\end{thm}

\begin{proof}
Given $f\in X_0+X_1$ and $u>0$, we consider the (quasi-) norms $(P_1 f)(u)$, $(R_1 f)(u)$ and $(Q_1 f)(u)$ defined as in \eqref{epqr} and we redefine 
\begin{align*}
(P_0 f)(u) &= \Big\| \b_0(t) \|s^{- \theta_0} \a_0(s) K(s,f)\|_{ \widetilde{F}_0(0,t)}\Big\|_{ \widetilde{E}_0(0,u)}, \\
(R_0 f)(u) &= \|\b_0\|_{\widetilde{E}_0(u,\infty) }\|t^{- \theta_0} \a_0(t) K(t,f)\|_{ \widetilde{F}_0(0,u)},\\
(Q_0 f)(u)&=\Big\| \b_0(t) \|s^{- \theta_0} \a_0(s) K(s,f)\|_{ \widetilde{F}_0(u,t)}\Big\|_{ \widetilde{E}_0(u,\infty)}.
\end{align*}
We denote, as usual, $Y_0= \overline{X}^{\mathcal L}_{\theta_0, \b_0,E_0,\a_0,F_0}$,  $Y_1= \overline{X}^{\mathcal R}_{\theta_1, \b_1,E_1,\a_1,F_1}$. We want to show that
\begin{equation}\label{5.1733}
K(\rho(u), f;Y_0,Y_1) \sim (P_{0}f)(u)+(R_0f)(u)+ \rho(u) [(R_1f)(u)+(Q_{1}f)(u)]
\end{equation}
where $\rho$ is defined by \eqref{rhoLR}. Again we follow the same steps as in Theorem \ref{5.12}.
We fix $f \in Y_0 + Y_1$ and  $u>0$.  We choose a decomposition $f = g + h$
such that $ \|g\|_{X_0} + u\|h\|_{X_1} \leq 2 K(u,f)$ and \eqref{4552} is satisfied.
Hence, to obtain the upper estimate of \eqref{5.1733} it is enough to prove that
$$
\|g\|_{Y_0}+\rho(u)\|h\|_{Y_1}\lesssim (P_0f)(u)+(R_0f)(u)+\rho(u) [(R_1f)(u)+(Q_{1}f)(u)].
$$
The inequality $\|h\|_{Y_1}\lesssim (R_1f)(u)+(Q_{1}f)(u)$ can be proved exactly as we did through \eqref{463} to \eqref{466}. Now, we are going to establish that $\|g\|_{Y_0}\lesssim (P_0 f)(u)+(R_0f)(u)$. By the triangle inequality and the (quasi-) subadditivity of the $K$-functional, we have that
\begin{align*}
\|g\|_{Y_0}&\leq(P_0g)(u)+(R_0g)(u)+(Q_0g)(u)\\&\lesssim (P_0f)(u)+(P_0h)(u)+(R_0f)(u)+(R_0h)(u)+(Q_0g)(u).
\end{align*}
For the proof of $(P_0 h)(u)\lesssim (P_0 f)(u)$ proceed as in \eqref{457}.
Besides that \eqref{4552}, Lemma \ref{lem1} (i) and \eqref{e1} give 
\begin{align*}
(R_0 h)(u) 
&\lesssim \frac{K(u,f)}{u} \| \b_0\|_{ \widetilde{E}_0(u,\infty)} \|t^{1- \theta_0} \a_0(t) \|_{ \widetilde{F}_0(0,u)}\nonumber\\
&\sim u^{-\theta_0}\a_0(u)\|\b_0\|_{\widetilde{E}_0(u,\infty)}K(u,f)\lesssim (R_0 f)(u).\nonumber
\end{align*}

Similarly, it follows
\begin{align*}
(Q_0 g)(u) 
&\lesssim K(u,f)\Big\| \b_0(t) \|s^{- \theta_0} \a_0(s)\|_{ \widetilde{F}_0(u,\infty)}\Big\|_{ \widetilde{E}_0(u,\infty)}\\&\sim u^{-\theta_0}\a_0(u)\|\b_0\|_{\widetilde{E}_0(u,\infty)}K(u,f)\lesssim (R_0 f)(u).
\end{align*}
Hence $\|g\|_{Y_0}\lesssim (P_0 f)(u)+(R_0 f)(u)$ and summing up
$$\|g\|_{Y_0}+\rho(u)\|h\|_{Y_1}\lesssim (P_0f)(u)+(R_0 f)(u)+\rho(u)[(R_1 f)(u)+(Q_{1}f)(u)].$$
This concludes the  proof of the upper estimate of (\ref{5.1733}).

Next, we proceed to the lower estimate of \eqref{5.1733}, that is 
$$(P_0f)(u)+(R_0 f)(u)+\rho(u)[(R_1 f)(u)+(Q_1f)(u)]\lesssim K(\rho(u), f;Y_0,Y_1)$$
for all $f\in Y_0+Y_1$ and $u>0$ where $\rho$ is defined by \eqref{rhoLR}.

We fix again $u>0$, $f\in Y_0+Y_1$. Let $f=f_0+f_1$ be any decomposition of $f$  with $f_0\in Y_0$ and $f_1\in Y_1$. Using the (quasi-) subadditivity of the $K$-functional and the definition of the norm in $Y_0$, $Y_1$, we have
\begin{align*}
(P_{0}f)(u) & \lesssim  (P_{0}f_0)(u) + (P_{0}f_1)(u)  \leq  \|f_0\|_{Y_{0}} + (P_{0}f_1)(u), \\
(R_{0}f)(u) & \lesssim  (R_{0}f_0)(u) + (R_{0}f_1)(u)  \leq  \|f_0\|_{Y_{0}} + (R_{0}f_1)(u),\\
(R_{1}f)(u) & \lesssim (R_{1}f_0)(u) + (R_{1}f_1)(u)  \leq(R_{1}f_0)(u) + \|f_1\|_{Y_{1}},\\
(Q_{1}f)(u) & \lesssim (Q_{1}f_0)(u) + (Q_{1}f_1)(u)  \leq
(Q_{1}f_0)(u) + \|f_1\|_{Y_{1}}.
\end{align*}
Thus, it is enough to verify that
\begin{equation}\label{870}
(P_{0}f_1)(u) + (R_0 f_1)(u)+\rho(u) [(R_1 f_0)(u)+(Q_{1}f_0)(u)]\lesssim \| f_0 \|_{Y_{0}} + \rho(u)\|f_1\|_{Y_{1}}.
\end{equation}
Arguing  as in \eqref{462}-\eqref{464bis} we can obtain that $(P_0f_1)(u)\lesssim \rho(u)\|f_1\|_{Y_1}$, $(R_1f_0)(u)\lesssim\frac1{\rho(u)}\|f_0\|_{Y_0}$ and $(Q_1f_0)(u)\lesssim\frac1{\rho(u)}\|f_0\|_{Y_0}$. On the other hand, using \eqref{eK2} with $f=f_1$ and Lemma \ref{lem1} (i), we have
\begin{align*}
(R_0f_1)(u)
&\lesssim \|f_1\|_{Y_1}\| \b_0\|_{\widetilde{E}_0(u,\infty)}\Big\|t^{\theta_1- \theta_0} \frac{\a_0(t)}{\a_1(t)\|\b_1\|_{\widetilde{E}_1(0,t)}}\Big\|_{ \widetilde{F}_0(0,u)}\\
&\sim u^{\theta_1- \theta_0} \frac{\a_0(u)\| \b_0\|_{\widetilde{E}_0(u,\infty)}}{\a_1(u)\|\b_1\|_{\widetilde{E}_1(0,u)}} \|f_1\|_{Y_1}=\rho(u)\|f_1\|_{Y_1}.
\end{align*}
Putting together the previous equations  we establish \eqref{870}.
Finally, by taking infimum over all possible decomposition of $f=f_0+f_1$, with $f_0\in Y_0$ and $f_1\in Y_1$, we obtain the desired estimate.
\end{proof}



\section{Reiteration formulae}\label{reiterationtheorems}
Our objective in this section is to identify the spaces 
$$(\overline{X}^{\mathcal R}_{\theta_0,\b_0,E_0,\a_0,F_0}, \overline{X}^{\mathcal R}_{\theta_1, \b_1,E_1,\a_1,F_1})_{\theta,\b,E},\quad
(\overline{X}^{\mathcal L}_{\theta_0, \b_0,E_0,\a_0,F_0}, \overline{X}^{\mathcal L}_{\theta_1,\b_1,E_1,\a_1,F_1})_{\theta,\b,E},$$
$$(\overline{X}^{\mathcal R}_{\theta_0,\b_0,E_0,\a_0,F_0}, \overline{X}^{\mathcal L}_{\theta_1, \b_1,E_1,\a_1,F_1})_{\theta,\b,E}\mand
(\overline{X}^{\mathcal L}_{\theta_0, \b_0,E_0,\a_0,F_0}, \overline{X}^{\mathcal R}_{\theta_1,\b_1,E_1,\a_1,F_1})_{\theta,\b,E}$$
for all possible values of $\theta\in[0,1]$. In that process we shall need two lemmas that we collect in next subsection.
\label{sereiteration}
\subsection{Lemmas}\hspace{1mm}

\vspace{2mm}
The first lemma is a change of variables and the second one is an equivalence between norms.
 
\begin{lem}\label{Le51}\cite[Lemma 4.1]{FMS-RL1}
Let $E$ be an r.i. space, $0<\alpha<1$, $\a$, $\b\in SV$ and consider the function $\rho(u)=u^{\alpha}\a(u)$, $u>0$. Then
$$\big\|  \rho(u)^{- \theta}  \b(\rho(u)) K( \rho(u), f)  \big\|_{\widetilde{E}}\sim\big\|  u^{- \theta}  \b(u) K( u, f)  \big\|_{\widetilde{E}}$$
for all $0\leq\theta\leq1$ and $f\in X_0+X_1$, with equivalent constant independent of $f$.
\end{lem}

\begin{lem}\label{lema45}
Let $E$, $F$ be r.i. spaces, $\a$, $\b\in SV$ and $\alpha, \beta, \gamma \in \R$ with
$\beta<0$ and $\gamma>0$. Then, the  equivalences
\begin{equation}\label{0t}
\Big \|  t^{\beta} \b(t)  \| s^{\alpha}  \a(s)  f(s)\|_{\widetilde{F}(0,t)} \Big \|_{\widetilde{E}}  \sim  \|  t^{\alpha +\beta} \a(t)  \b(t)  f(t) \|_{\widetilde{E}}
\end{equation}
and
\begin{equation}\label{tinfty}
\Big \|  t^{\gamma} \b(t)  \| s^{\alpha}  \a(s)  f(s)\|_{\widetilde{F}(t, \infty)} \Big \|_{\widetilde{E}}  \sim  \|  t^{\alpha + \gamma} \a(t)  \b(t)  f(t) \|_{\widetilde{E}}
\end{equation}
hold for any monotone measurable function $f$ on $(0,\infty)$.
\end{lem}

\begin{proof}
Both properties are proved in \cite[Theorems 3.6 and 3.7]{FMS-4} for non-increasing measurable functions on $(0,\infty)$. We assume now that $f$ is a non-decreasing measurable function. By \eqref{fbar} and applying \eqref{tinfty}, with $\overline{f}$ and $\gamma=-\beta<0$, we obtain \eqref{0t} for $f$
\begin{align*}
\Big \|  t^{\beta} \b(t)  \| s^{\alpha}  \a(s)  f(s)\|_{\widetilde{F}(0,t)} \Big \|_{\widetilde{E}}  &=
 \Big \|  t^{\beta} \b(t)  \| s^{-\alpha}  \overline{\a}(s)  \overline{f}(s)\|_{\widetilde{F}(1/t,\infty)} \Big \|_{\widetilde{E}}\\
&=
 \Big \|  t^{-\beta} \overline{\b}(t)  \| s^{-\alpha}  \overline{\a}(s)  \overline{f}(s)\|_{\widetilde{F}(t,\infty)} \Big \|_{\widetilde{E}}\\
&\sim \|  t^{-(\alpha + \beta)} \overline{\a}(t)  \overline{\b}(t)  \overline{f}(t) \|_{\widetilde{E}}\\
&=  \|  t^{\alpha + \beta} \a(t)  \b(t)  f(t) \|_{\widetilde{E}}.
\end{align*}
The proof of \eqref{tinfty} can be done in a similar way.
\end{proof}

Now, we are in position to establish reiteration theorems for ${\mathcal R}$ and ${\mathcal L}$-classes.

\subsection{The space $(\overline{X}^{\mathcal R}_{\theta_0, \b_0,E_0,\a_0,F_0}, \overline{X}^{\mathcal R}_{\theta_1,\b_1,E_1,\a_1,F_1})_{\theta,\b,E}$, $0<\theta_0<\theta_1<1$ and $0\leq \theta\leq 1$.} 

\begin{thm} \label{thmRR}
Let $0<\theta_0<\theta_1<1$. Let $E$, $E_0$, $E_1$, $F_0$, $F_1$ r.i. spaces, $\a_0$, $\a_1$, $\b$,  $\b_0$, $\b_1\in SV$   such that $\|\b_0\|_{\widetilde{E}_0(0,1)}<\infty$ and  $\|\b_1\|_{\widetilde{E}_1(0,1)}<\infty$ and denote
$$\rho(u) = u^{\theta_1- \theta_0} \frac{\a_0(u)\|\b_0\|_{\widetilde{E}_0(0,u)}}{\a_1(u)\|\b_1\|_{\widetilde{E}_1(0,u)} },\quad u>0.$$
\begin{itemize}
\item[a)] If $0<\theta<1$, then
\begin{equation*}\label{reRR}
(\overline{X}^{\mathcal R}_{\theta_0, \b_0,E_0,\a_0,F_0}, \overline{X}^{\mathcal R}_{\theta_1,\b_1,E_1,\a_1,F_1})_{\theta,\b,E}= \overline{X}_{\tilde{\theta},\B_\theta,E}
\end{equation*} 
where $\tilde{\theta}=(1-\theta)\theta_0+\theta\theta_1$
and
$$B_\theta(u)=(\a_0(u)\|\b_0\|_{\widetilde{E}_0(0,u)})^{1-\theta}(\a_1(u)\|\b_1\|_{\widetilde{E}_1(0,u)})^{\theta}\b(\rho(u)),\quad u>0.$$ 

\item[b)] If $\theta=0$ and $\|\b\|_{\widetilde{E}(1,\infty)}<\infty$, then
$$(\overline{X}^{\mathcal R}_{\theta_0, \b_0,E_0,\a_0,F_0}, \overline{X}^{\mathcal R}_{\theta_1,\b_1,E_1,\a_1,F_1})_{0,\b,E}= \overline{X}^{\mathcal R,\mathcal L}_{\theta_0,\b\circ\rho,E,\b_0,E_0,a_0,F_0}.$$

\item[c)]  If $\theta=1$ and $\|\b\|_{\widetilde{E}(0,1)}<\infty$, then 
$$(\overline{X}^{\mathcal R}_{\theta_0, \b_0,E_0,\a_0,F_0}, \overline{X}^{\mathcal R}_{\theta_1,\b_1,E_1,\a_1,F_1})_{1,\b,E}=\overline{X}^{\mathcal R}_{\theta_1,\B_1,E,\a_1,F_1}\cap\overline{X}^{\mathcal R,\mathcal R}_{\theta_1,\b\circ\rho,E,\b_1,E_1,a_1,F_1}$$
\vspace{1mm}
where $\B_1(u)=\|\b_1\|_{\widetilde{E}_1(0,u)}\b(\rho(u))$, $u>0$.
\end{itemize}
\end{thm}

\begin{proof}
Throughout the proof we use the notation $Y_0=\overline{X}^{\mathcal R}_{\theta_0, \b_0,E_0,\a_0,F_0}$, $Y_1=\overline{X}^{\mathcal R}_{\theta_1, \b_1, E_1,\a_1,F_1}$ and $\overline{K}(u,f)= K(u, f;Y_0,Y_1)$, $u>0$.

We start with the proof of a). Let  $f \in \overline{Y}_{\theta,\b,E}$. The generalized Holmstedt type formula \eqref{eqHoRR} and estimate \eqref{e5} give that
\begin{align*}
\overline{K}(\rho(u),f)&
 \gtrsim\Big\| \b_0(t) \|s^{- \theta_0} \a_0(s) K(s,f)\|_{ \widetilde{F}_0(t,u)}\Big\|_{ \widetilde{E}_0(0,u)}\label{900}\\ &\gtrsim u^{-\theta_0}\a_0(u)\|\b_0\|_{\widetilde{E}_0(0,u)}K(u,f).\nonumber
\end{align*}
Hence, using also Lemma \ref{Le51}, we have
\begin{align*}
\|f\|_{\overline{Y}_{\theta,\b,E}}&\sim \|\rho(u)^{-\theta}\b(\rho(u))\overline{K}(\rho(u), f)\|_{\widetilde{E}}\\
&\gtrsim \big\|\rho(u)^{-\theta}\b(\rho(u))u^{- \theta_0} \a_0(u)\|\b_0\|_{\widetilde{E}_0(0,u)} K(u,f) \big \|_{\widetilde{E}}.
\end{align*}
Observing  
\begin{equation}\label{rhob}
\rho(u)^{- \theta}  \b(\rho(u))=u^{\theta_0-\tilde{\theta}}\frac{\B_\theta(u)}{ \a_0(u)\|\b_0\|_{\widetilde{E}_0(0,u)}},\quad u>0,
\end{equation}
we obtain that $\|f\|_{\overline{X}_{\tilde{\theta},\B_\theta,E}}\lesssim \|f\|_{\overline{Y}_{\theta,\b,E}}$ . This proves the the inclusion $\overline{Y}_{\theta,\b,E}\hookrightarrow\overline{X}_{\tilde{\theta},\B_\theta,E}$.

Next, we proceed with the reverse inclusion. Let $f\in \overline{X}_{\tilde{\theta},\B_\theta,E}$. Using again Lemma \ref{Le51}, the generalized Holmstedt type formula \eqref{eqHoRR} and the triangular inequality, we have
\begin{align*}
\|f\|_{\overline{Y}_{\theta,\b,E}}&\sim\big\|  \rho(u)^{-\theta} \b(\rho(u)) \overline{K}(\rho(u),f)
\big\|_{\widetilde{E}} \\
 & \lesssim 
 \Big \|\rho(u)^{-\theta} \b(\rho(u))\big\| \b_0(t) 
     \|s^{-\theta_0}\a_0(s)K(s,f)\|_{\widetilde{F}_0(t,u)}\big\|_{\widetilde{E}_{0}(0,u)}  \Big\|_{\widetilde{E}}\\
 &+  \Big\|\rho(u)^{1-\theta} \b(\rho(u)) \| \b_1\|_{\widetilde{E}_1(0,u)}
     \|t^{-\theta_1}\a_1(t)K(t,f)\|_{\widetilde{F}_1(u,\infty)}\Big\|_{\widetilde{E}} 
 \\&+\Big\|\rho(u)^{1-\theta} \b(\rho(u))  \big\| \b_1(t) 
     \|s^{-\theta_1}\a_1(s)K(s,f)\|_{\widetilde{F}_1(t,\infty)}\big\|_{\widetilde{E}_{1}(u,\infty)}   \Big\|_{\widetilde{E}}\\
		&:=I_1+I_2+I_3. \end{align*}
Hence, in order to prove that $f\in\overline{Y}_{\theta,\b,E}$ it is suffices to estimate last three expressions, $I_1$, $I_2$ and $I_3$, by $\| u^{- \tilde{\theta}} \B_\theta(u) K(u,f) \|_{\widetilde{E}}$. We start with $I_1$. 
Using \eqref{rhob} and Lemma \ref{lema45}, with the $K$-functional and  $\beta=\theta_0-\tilde{\theta}<0$, we obtain
 \begin{align}
 I_1&=\bigg \|u^{\theta_{0} - \tilde{\theta}} \tfrac{\B_\theta(u)}{\a_0(u)\|\b_0\|_{\widetilde{E}_0(0,u)}}\Big\| \b_0(t) 
     \|s^{-\theta_0}\a_0(s)K(s,f)\|_{\widetilde{F}_0(t,u)}\Big\|_{\widetilde{E}_{0}(0,u)}     \bigg \|_{\widetilde{E}}\label{I1}\\
     &\leq \bigg \|  u^{\theta_{0} - \tilde{\theta}} \tfrac{\B_\theta(u)}{\a_0(u)} \|s^{-\theta_0}\a_0(s)K(s,f)\|_{\widetilde{F}_0(0,u)}  \bigg \|_{\widetilde{E}}\nonumber\\
      &\sim\big \| u^{- \tilde{\theta}} \B_\theta(u) K(u,f)\big \|_{\widetilde{E}}.\nonumber
\end{align}
Now we estimate $I_2$. The relation
\begin{equation}\label{erho2}
\rho(u)^{1- \theta}  \b(\rho(u))=u^{\theta_1-\tilde{\theta}}\frac{\B_\theta(u)}{ \a_1(u)\|\b_1\|_{\widetilde{E}_1(0,u)}},\quad u>0,\end{equation}
and Lemma \ref{lema45}, with the $K$-functional and  $\gamma=\theta_1-\widetilde{\theta}>0$, give that 
\begin{equation*}\label{I2}
 I_2=\Big \|u^{\theta_1-\tilde{\theta}}\tfrac{\B_\theta(u)}{\a_1(u)}
     \|t^{-\theta_1}\a_1(t)K(t,f)\|_{\widetilde{F}_1(u,\infty)}  \Big \|_{\widetilde{E}}
  \sim \|u^{-\tilde{\theta}}\B_\theta(u)K(u,f\|_{\widetilde{E}}.
\end{equation*}
Finally, we approach $I_3$ through $I_2$.  Using \eqref{erho2} we have
$$I_3=\bigg \|u^{\theta_1-\tilde{\theta}}\tfrac{\B_\theta(u)}{\a_1(u)\|\b_1\|_{\widetilde{E}_1(0,u)}}\big\| \b_1(t) \|s^{-\theta_1}\a_1(s)K(s,f)\|_{\widetilde{F}_1(t,\infty)}\big\|_{\widetilde{E}_{1}(u,\infty)} \bigg \|_{\widetilde{E}}.$$
Since the function $t\rightsquigarrow\|\cdot\|_{\widetilde{F}_1(t,\infty)}$ is monotone and $\theta_1-\tilde{\theta}>0$, we can apply Lemma \ref{lema45} to obtain that 
$$
I_3 \sim \bigg \|u^{\theta_1-\tilde{\theta}}\tfrac{\B_\theta(u)}{\a_1(u)\|\b_1\|_{\widetilde{E}_1(0,u)}} \b_1(u) \|s^{-\theta_1}\a_1(s)K(s,f)\|_{\widetilde{F}_1(u,\infty)} \bigg \|_{\widetilde{E}}.$$
From Lemma \ref{lem1} (iv) and \eqref{erho2} we deduce that $I_3\lesssim I_2$ and the proof of a) is completed.

The proof of b) follows similar steps. The inclusion $\overline{Y}_{0,\b,E}\hookrightarrow \overline{X}^{\mathcal R,\mathcal L}_{\theta_0,\b\circ\rho,E,\b_0,E_0,a_0,F_0}$ comes directly from Theorem \ref{5.12} and the lattice property of $\widetilde{E}$.
In order to prove the reverse inclusion we use again Lemma \ref{Le51}, Theorem \ref{5.12} and the triangular inequality. Indeed, let $f\in\overline{X}^{\mathcal R,\mathcal L}_{\theta_0,\b\circ\rho,E,\b_0,E_0,a_0,F_0}$, then 
\begin{align*}
\|f\|_{\overline{Y}_{0,\b,E}}&\sim\big\|  \b(\rho(u))\overline{K}(\rho(u),f)
\big\|_{\widetilde{E}} \\
 & \lesssim 
 \Big \|\b(\rho(u))\big\| \b_0(t) 
     \|s^{-\theta_0}\a_0(s)K(s,f)\|_{\widetilde{F}_0(t,u)}\big\|_{\widetilde{E}_{0}(0,u)}  \Big\|_{\widetilde{E}}\\
 &+  \Big\|\rho(u) \b(\rho(u)) \| \b_1\|_{\widetilde{E}_1(0,u)}
     \|t^{-\theta_1}\a_1(t)K(t,f)\|_{\widetilde{F}_1(u,\infty)}\Big\|_{\widetilde{E}} 
 \\&+\Big\|\rho(u)\b(\rho(u))  \big\| \b_1(t) 
     \|s^{-\theta_1}\a_1(s)K(s,f)\|_{\widetilde{F}_1(t,\infty)}\big\|_{\widetilde{E}_{1}(u,\infty)} \Big\|_{\widetilde{E}}\\&:=I_4+I_5+I_6. \end{align*}
It is clear that $I_4=\|f\|_{\overline{X}^{\mathcal R,\mathcal L}_{\theta_0,\b\circ\rho,E,\b_0,E_0,a_0,F_0}}$. Let us estimate $I_5$ by $I_4$. Lemma \ref{lema45}, with the $K$-functional and  $\gamma=\theta_1-\theta_0>0$,
 and \eqref{e5} guarantee that
\begin{align*}
I_5&= \Big\|u^{\theta_1-\theta_0}\frac{\a_0(u)\| \b_0\|_{\widetilde{E}_0(0,u)}}{\a_1(u)} \b(\rho(u)) 
     \|t^{-\theta_1}\a_1(t)K(t,f)\|_{\widetilde{F}_1(u,\infty)}\Big\|_{\widetilde{E}} \\
		&\sim\Big\|u^{-\theta_0}\a_0(u)\| \b_0\|_{\widetilde{E}_0(0,u)} \b(\rho(u)) 
     K(u,f)\Big\|_{\widetilde{E}}\lesssim I_4.
\end{align*}
The third term $I_6$ is bounded by $I_5$. Indeed, since the function $t\rightsquigarrow \|\cdot\|_{ \widetilde{F}_1(t,\infty)}$ is monotone and $\gamma=\theta_1-\tilde{\theta}>0$, Lemma \ref{lema45}  gives that
\begin{align*}
I_6&=\Bigg\|u^{\theta_1-\theta_0}\frac{\a_0(u)\|\b_0\|_{\widetilde{E}_0(0,u)}}{\a_1(u)\|\b_1\|_{\widetilde{E}_1(0,u)}}\b(\rho(u))  \Big\| \b_1(t) 
     \|s^{-\theta_1}\a_1(s)K(s,f)\|_{\widetilde{F}_1(t,\infty)}\Big\|_{\widetilde{E}_{1}(u,\infty)} \Bigg\|_{\widetilde{E}}\\
		&\sim \Bigg\|u^{\theta_1-\theta_0}\frac{\a_0(u)\|\b_0\|_{\widetilde{E}_0(0,u)}}{\a_1(u)\|\b_1\|_{\widetilde{E}_1(0,u)}}\b(\rho(u))  \b_1(u) 
     \|s^{-\theta_1}\a_1(s)K(s,f)\|_{\widetilde{F}_1(u,\infty)}\Bigg\|_{\widetilde{E}}.
\end{align*}
Then, Lemma \ref{lem1} (iv) implies that $I_6\lesssim I_5$. Summing up, $$\|f\|_{\overline{Y}_{0,\b,E}}\lesssim I_4+I_5+I_6\lesssim I_4=\|f\|_{\overline{X}^{\mathcal R,\mathcal L}_{\theta_0,\b\circ\rho,E,\b_0,E_0,a_0,F_0}}$$ and the proof of b) is concluded.

Finally, we proceed with the proof of c).  Let $f\in \overline{Y}_{1,\b,E}$. Lemma \ref{Le51},  
Theorem \ref{5.12} and the lattice property guarantee that 
\begin{align*}
\|f\|_{\overline{Y}_{1,\b,E}}&\sim
\|\rho(u)^{-1}\b(\rho(u))\overline{K}(\rho(u), f)\|_{\widetilde{E}}\\
&\gtrsim \Big\|\b(\rho(u))\, \| \b_1\|_{\widetilde{E}_1(0,u)} \|t^{- \theta_1} \a(t) K(t,f)\|_{ \widetilde{F}_1(u,\infty)}\Big\|_{\widetilde{E}}=\|f\|_{\overline{X}^{\mathcal R}_{\theta_1,\B_1,E,\a_1,F_1}}
\end{align*}
and 
$$
\|f\|_{\overline{Y}_{1,\b,E}}\gtrsim\bigg\|\b(\rho(u))\,\Big\| \b_1(t) \|s^{- \theta_1} \a(s) K(s,f)\|_{ \widetilde{F}_1(t,\infty)}\Big\|_{ \widetilde{E}_1(u,\infty)}\bigg\|_{\widetilde{E}}=\|f\|_{\overline{X}^{\mathcal R,\mathcal R}_{\theta_1,\b\circ\rho,E,\b_1,E_1,a_1,F_1}}.$$
Then, $f\in \overline{X}^{\mathcal R}_{\theta_1,\B_1,E,\a_1,F_1}\cap \overline{X}^{\mathcal R,\mathcal R}_{\theta_1,\b\circ\rho,E,\b_1,E_1,a_1,F_1}$.

Next, we study the reverse inclusion. Let $f\in \overline{X}^{\mathcal R}_{\theta_1,\B_1,E,\a_1,F_1}\cap \overline{X}^{\mathcal R,\mathcal R}_{\theta_1,\b\circ\rho,E,\b_1,E_1,a_1,F_1}$. Using again Lemma \ref{Le51}, Theorem \ref{5.12} and the triangular inequality, we have that
\begin{align*}
\|f\|_{\overline{Y}_{1,\b,E}}&\sim\|\rho(u)^{-1}\b(\rho(u))\overline{K}(\rho(u), f)\|_{\widetilde{E}}
\\&\lesssim
 \Big \|\rho(u)^{-1} \b(\rho(u))\big\| \b_0(t) 
     \|s^{-\theta_0}\a_0(s)K(s,f)\|_{\widetilde{F}_0(t,u)}\big\|_{\widetilde{E}_{0}(0,u)}  \Big\|_{\widetilde{E}}\\
 &+  \Big\|\b(\rho(u)) \| \b_1\|_{\widetilde{E}_1(0,u)}
     \|t^{-\theta_1}\a_1(t)K(t,f)\|_{\widetilde{F}_1(u,\infty)}\Big\|_{\widetilde{E}} 
 \\&+\Big\|\b(\rho(u))  \big\| \b_1(t) 
     \|s^{-\theta_1}\a_1(s)K(s,f)\|_{\widetilde{F}_1(t,\infty)}\big\|_{\widetilde{E}_{1}(u,\infty)}   \Bigg\|_{\widetilde{E}}\\
		&:=I_7+I_8+I_9.
\end{align*}
Since $I_8=\|f\|_{\overline{X}^{\mathcal R}_{\theta_1,B_1,E,a_1,F_1}}$ and $I_9=\|f\|_{\overline{X}^{\mathcal R,\mathcal R}_{\theta_1,\b\circ\rho,E,\b_1,E_1,a_1,F_1}}$, it is enough to estimate $I_7$. Using Lemma \ref{lema45}, with the $K$-functional and $\beta=\theta_0-\theta_1<0$, and estimate \eqref{e3}, we obtain 
\begin{align*}
I_7&=\bigg\|u^{\theta_0-\theta_1}\frac{\a_1(u)\|\b_1\|_{\widetilde{E}_1(0,u)}}{\a_0(u)\|\b_0\|_{\widetilde{E}_0(0,u)}}\,\b(\rho(u)) \big\| \b_0(t) 
     \|s^{-\theta_0}\a_0(s)K(s,f)\|_{\widetilde{F}_0(t,u)}\big\|_{\widetilde{E}_{0}(0,u)}\bigg\|_{\widetilde{E}}\\
&\leq	\bigg\|u^{\theta_0-\theta_1}\frac{\a_1(u)\|\b_1\|_{\widetilde{E}_1(0,u)}}{\a_0(u)}\,\b(\rho(u))
     \|s^{-\theta_0}\a_0(s)K(s,f)\|_{\widetilde{F}_0(0,u)}\bigg\|_{\widetilde{E}}\\	
&\sim	\bigg\|u^{-\theta_1}\a_1(u)\|\b_1\|_{\widetilde{E}_1(0,u)}\,\b(\rho(u))
     K(u,f)\bigg\|_{\widetilde{E}}\\
&\lesssim	\bigg\|\b(\rho(u))\|\b_1\|_{\widetilde{E}_1(0,u)}
    \|t^{-\theta_1}\a_1(t)K(t,f)\|_{\widetilde{F}_1(u,\infty)} \bigg\|_{\widetilde{E}}=I_8=\|f\|_{\overline{X}^{\mathcal R}_{\theta_1,B_1,E,a_1,F_1}}.
\end{align*}		
This completes the proof.
\end{proof} 

\subsection{The space $(\overline{X}^{\mathcal L}_{\theta_0, \b_0,E_0,a_0,F_0}, \overline{X}^{\mathcal L}_{\theta_1,\b_1,E_1,a_1,F_1})_{\theta,\b,E}$, $0<\theta_0<\theta_1<1$ and $0\leq \theta\leq 1$.}\hspace{1mm}

\vspace{2mm}
Next theorem can be proved similarly, except that Holmstedt formula \eqref{eHL} has to be used in place of \eqref{eqHoRR}. However, since this approach is  lenghty, we prove it using Theorem \ref{thmRR} and a symmetry argument.

\begin{thm} \label{thmLL}
Let $0<\theta_0<\theta_1<1$. Let $E$, $E_0$, $E_1$, $F_0$, $F_1$ r.i. spaces, $a_0$, $a_1$, $\b$, $\b_0$, $\b_1\in SV$ such that $\|\b_0\|_{\widetilde{E}_0(1,\infty)}<\infty$ and  $\|\b_1\|_{\widetilde{E}_1(1,\infty)}<\infty$ and denote
$$\rho(u) = u^{\theta_1- \theta_0} \frac{\a_0(u)\|\b_0\|_{\widetilde{E}_0(u,\infty)}}{a_1(u)\|\b_1\|_{\widetilde{E}_1(u,\infty)} },\quad u>0.$$
\begin{itemize}
\item[a)] If $0<\theta<1$, then
$$
(\overline{X}^{\mathcal L}_{\theta_0, \b_0,E_0,a_0,F_0}, \overline{X}^{\mathcal L}_{\theta_1,\b_1,E_1,a_1,F_1})_{\theta,\b,E}= \overline{X}_{\tilde{\theta},\B_\theta,E}
$$
where $\tilde{\theta}=(1-\theta)\theta_0+\theta\theta_1$
and
$$B_\theta(u)=(a_0(u)\|\b_0\|_{\widetilde{E}_0(u,\infty)})^{1-\theta}(a_1(u)\|\b_1\|_{\widetilde{E}_1(u,\infty)})^{\theta}\b(\rho(u)),\quad u>0.$$ 

\item[b)] If $\theta=0$ and $\|\b\|_{\widetilde{E}(1,\infty)}<\infty$, then
$$(\overline{X}^{\mathcal L}_{\theta_0, \b_0,E_0,a_0,F_0}, \overline{X}^{\mathcal L}_{\theta_1,\b_1,E_1,a_1,F_1})_{0,\b,E}=\overline{X}^{\mathcal L}_{\theta_0,\B_0,E,a_0,F_0}\cap\overline{X}^{\mathcal L,\mathcal L}_{\theta_0,\b\circ\rho,E,\b_0,E_0,a_0,F_0}$$
\vspace{1mm}
where $\B_0(u)=\|\b_0\|_{\widetilde{E}_0(u,\infty)}\b(\rho(u))$, $u>0$.
\item[c)]  If $\theta=1$ and $\|\b\|_{\widetilde{E}(0,1)}<\infty$, then 
\begin{equation}\label{theta1}
(\overline{X}^{\mathcal L}_{\theta_0, \b_0,E_0,a_0,F_0}, \overline{X}^{\mathcal L}_{\theta_1,\b_1,E_1,a_1,F_1})_{1,\b,E}=
\overline{X}^{\mathcal L,\mathcal R}_{\theta_1,\b\circ\rho,E,\b_1,E_1,a_1,F_1}.
\end{equation}
\end{itemize}
\end{thm}

\begin{proof} We prove c). 
Applying both identities of Lemma \ref{symLR} we have
\begin{align*}
(\overline{X}^{\mathcal L}_{\theta_0, \b_0,E_0,a_0,F_0}, \overline{X}^{\mathcal L}_{\theta_1,\b_1,E_1,\a_1,F_1})_{1,\b,E}
&=(\overline{X}^{\mathcal L}_{\theta_1,\b_1,E_1,a_1,F_1},\overline{X}^{\mathcal L}_{\theta_0, \b_0,E_0,a_0,F_0})_{0,\overline{\b},E}\\
&=((X_1,X_0)^{\mathcal R}_{1-\theta_1,\overline{\b}_1,E_1,\overline{a}_1,F_1},
(X_1,X_0)^{\mathcal R}_{1-\theta_0,\overline{\b}_0,E_0,\overline{a}_0,F_0})_{0,\overline{\b},E}.
\end{align*} 
Then, by Theorem \ref{thmRR} b) we obtain
$$(\overline{X}^{\mathcal L}_{\theta_0, \b_0,E_0,a_0,F_0}, \overline{X}^{\mathcal L}_{\theta_1,\b_1,E_1,a_1,F_1})_{1,\b,E}
=(X_1,X_0)^{\mathcal R,\mathcal L}_{1-\theta_1,\overline{\b}\circ\rho^{\#},E,\overline{\b}_1,E_1,\overline{a}_1,F_1}$$
where $$\rho^{\#}(u)=u^{\theta_1-\theta_0}\frac{\overline{a}_1(u)\|\overline{\b}_1\|_{\widetilde{E}_1(0,u)}}{\overline{a}_0(u)\|\overline{\b}_0\|_{\widetilde{E}_0(0,u)} },\quad u>0.$$
Equality\eqref{theta1} follows from Lemma \ref{symLL} and that $\overline{\overline{\b}\circ\rho^{\#}}(u)=\b\circ\rho(u)$, $u>0$.

The proofs of a) and b) can be carried out using similar arguments and Theorem \ref{thmRR} a), c).
\end{proof}

\subsection{The space $(\overline{X}^{\mathcal R}_{\theta_0, \b_0,E_0,\a_0,F_0}, \overline{X}^{\mathcal L}_{\theta_1,\b_1,E_1,a_1,F_1})_{\theta,\b,E}$, $0<\theta_0<\theta_1<1$ and $0\leq \theta\leq 1$.}\hspace{1mm}

\vspace{2mm}
The following theorems deal with reiteration formulae with mixed ${\mathcal R}$ and ${\mathcal L}$-spaces. Both theorems have to be proved independently, since the symmetry argument cannot be applied.

\begin{thm} \label{thmRL}
Let $0<\theta_0<\theta_1<1$. Let $E$, $E_0$, $E_1$, $F_0$, $F_1$ r.i. spaces, $a_0$, $a_1$, $\b$,  $\b_0$, $\b_1\in SV$   such that $\|\b_0\|_{\widetilde{E}_0(0,1)}<\infty$ and  $\|\b_1\|_{\widetilde{E}_1(1,\infty)}<\infty$ and denote
$$\rho(u) = u^{\theta_1- \theta_0} \frac{\a_0(u)\|\b_0\|_{\widetilde{E}_0(0,u)}}{\a_1(u)\|\b_1\|_{\widetilde{E}_1(u,\infty)}},\quad u>0.$$
\begin{itemize}
\item[a)] If $0<\theta<1$, then
$$
(\overline{X}^{\mathcal R}_{\theta_0, \b_0,E_0,a_0,F_0}, \overline{X}^{\mathcal L}_{\theta_1,\b_1,E_1,a_1,F_1})_{\theta,\b,E}= \overline{X}_{\tilde{\theta},\B_\theta,E}
$$
where $\tilde{\theta}=(1-\theta)\theta_0+\theta\theta_1$
and
$$B_\theta(u)=(a_0(u)\|\b_0\|_{\widetilde{E}_0(0,u)})^{1-\theta}(a_1(u)\|\b_1\|_{\widetilde{E}_1(u,\infty)})^{\theta}\b(\rho(u)),\quad u>0.$$ 

\item[b)] If $\theta=0$ and $\|\b\|_{\widetilde{E}(1,\infty)}<\infty$, then
$$(\overline{X}^{\mathcal R}_{\theta_0, \b_0,E_0,a_0,F_0}, \overline{X}^{\mathcal L}_{\theta_1,\b_1,E_1,a_1,F_1})_{0,\b,E}=\overline{X}^{\mathcal R,\mathcal L}_{\theta_0,\b\circ\rho,E,\b_0,E_0,a_0,F_0}.$$

\item[c)]  If $\theta=1$ and  $\|\b\|_{\widetilde{E}(0,1)}<\infty$, then 
$$(\overline{X}^{\mathcal R}_{\theta_0, \b_0,E_0,a_0,F_0}, \overline{X}^{\mathcal L}_{\theta_1,\b_1,E_1,a_1,F_1})_{1,\b,E}=
\overline{X}^{\mathcal L,\mathcal R}_{\theta_1,\b\circ\rho,E,\b_1,E_1,a_1,F_1}.$$
\end{itemize}
\end{thm}

\begin{proof}
The proof follows essentially the same steps as the proof of Theorem \ref{thmRR}, although some estimates are different.
Let us start by denoting as usual $Y_0=\overline{X}^{\mathcal R}_{\theta_0, \b_0,E_0,a_0,F_0}$, $Y_1=\overline{X}^{\mathcal L}_{\theta_1,\b_1,E_1,a_1,F_1}$ and $\overline{K}(u,f)=K(u,f;Y_0,Y_1)$, $u>0$.

The inclusion 
$\overline{Y}_{\theta,\b,E}\hookrightarrow \overline{X}_{\tilde{\theta},\B_\theta,E}$ follows from Lemma \ref{Le51}, Theorem \ref{5.17} and inequality \eqref{e5} as we did in the proof of Theorem \ref{thmRR} a). Let now $f\in \overline{X}_{\tilde{\theta},\B_\theta,E}$. Using Lemma \ref{Le51}, Theorem \ref{5.17} and the triangular inequality, we have that
\begin{align*}
\|f\|_{\overline{Y}_{\theta,\b,E}}&\sim\big\|  \rho(u)^{-\theta} \b(\rho(u)) \overline{K}(\rho(u),f)
\big\|_{\widetilde{E}} \\
 & \lesssim 
 \Big \|\rho(u)^{-\theta} \b(\rho(u))\big\| \b_0(t) 
     \|s^{-\theta_0}a_0(s)K(s,f)\|_{\widetilde{F}_0(t,u)}\big\|_{\widetilde{E}_{0}(0,u)}  \Big\|_{\widetilde{E}}\\
 &+  \Big\|\rho(u)^{1-\theta} \b(\rho(u))\big\| \b_1(t)\| s^{- \theta_1} \a_1(s) K(s,f)\|_{ \widetilde{F}_1(u,t)}\big\|_{\widetilde{E}_1(u,\infty)} \Big\|_{\widetilde{E}} 
\\
&:=I_{10}+I_{11}.
  \end{align*}
The expression $I_{10}$ can be estimate by $\|f\|_{\overline{X}_{\tilde{\theta},\B_\theta,E}}$ as we did with $I_1$ in \eqref{I1}. We also observe that
\begin{align*}
I_{11}&=\Bigg\|u^{\theta_1-\tilde{\theta}}\frac{B_\theta(u)}{a_1(u)\|\b_1\|_{\widetilde{E}_1(u,\infty)}}\Big\| \b_1(t)\| s^{- \theta_1} a_1(s) K(s,f)\|_{ \widetilde{F}_1(u,t)}\Big\|_{\widetilde{E}_1(u,\infty)} \Bigg\|_{\widetilde{E}}\label{902}\\
&\leq\Big\|u^{\theta_1-\tilde{\theta}}\frac{B_\theta(u)}{a_1(u)}\| s^{- \theta_1} \a_1(s) K(s,f)\|_{ \widetilde{F}_1(u,\infty)} \Big\|_{\widetilde{E}}.\nonumber
\end{align*}
Applying, as usual, Lemma \ref{lema45} with the $K$-functional and $\gamma=\theta_1-\widetilde{\theta}>0$, we  obtain that $I_{11}\lesssim \|f\|_{\overline{X}_{\tilde{\theta},\B_\theta,E}}$. This completes the proof of a). 

Next we proceed with the proof of b). The inclusion $\overline{Y}_{0,\b,E}\hookrightarrow \overline{X}^{\mathcal R,\mathcal L}_{\theta_0,\b\circ\rho,E,\b_0,E_0,a_0,F_0}$ follows directly from Lemma \ref{Le51}, Theorem \ref{5.17} and the lattice property of $\widetilde{E}$.
In order to prove the reverse inclusion we use again Lemma \ref{Le51}, Theorem \ref{5.17} and the triangular inequality to obtain that 
\begin{align*}
\|f\|_{\overline{Y}_{0,\b,E}}&\sim\big\|  \b(\rho(u))\overline{K}(\rho(u),f)
\big\|_{\widetilde{E}} \\
 & \lesssim 
 \Big \|\b(\rho(u))\big\| \b_0(t) 
     \|s^{-\theta_0}\a_0(s)K(s,f)\|_{\widetilde{F}_0(t,u)}\big\|_{\widetilde{E}_{0}(0,u)}  \Big\|_{\widetilde{E}}\\
&+\Big\|\rho(u)\b(\rho(u))  \big\| \b_1(t) 
     \|s^{-\theta_1}\a_1(s)K(s,f)\|_{\widetilde{F}_1(u,t)}\big\|_{\widetilde{E}_{1}(u,\infty)} \Big\|_{\widetilde{E}}\\&:=I_{12}+I_{13}. \end{align*}
It is clear that $I_{12}$ is equal to $\|f\|_{\overline{X}^{\mathcal R,\mathcal L}_{\theta_0,\b\circ\rho,E,\b_0,E_0,a_0,F_0}}$. Besides that, arguing as we did with $I_{11}$ and using \eqref{e5} we estimate $I_{13}$ by $I_{12}$. Indeed,
\begin{align*}
I_{13}
&\lesssim \big\|u^{-\theta_0}a_0(u)\|\b_0\|_{\widetilde{E}_0(0,u)}\b(\rho(u))K(u,f)\big\|_{\widetilde{E}}\\
&\lesssim \Big\|\b(\rho(u))\big\|\b_0(t)\|s^{-\theta_0} \a_0(s) K(s,f)\|_{\widetilde{F}_0(t,u)} \big\|_{\widetilde{E}_0(0,u)}\Big\|_{\widetilde{E}}=I_{12}.
\end{align*}
This completes the proof of b). Similar steps using \eqref{e7} instead of \eqref{e5} lead to the proof of c).
\end{proof}

\subsection{The space $(\overline{X}^{\mathcal L}_{\theta_0, \b_0,E_0,\a_0,F_0}, \overline{X}^{\mathcal R}_{\theta_1,\b_1,E_1,\a_1,F_1})_{\theta,\b,E}$, $0<\theta_0<\theta_1<1$ and $0\leq \theta\leq 1$.}

\begin{thm} \label{thmLR}
Let $0<\theta_0<\theta_1<1$. Let $E$, $E_0$, $E_1$, $F_0$, $F_1$ r.i. spaces, $\a_0$, $\a_1$, $\b$,  $\b_0$, $\b_1\in SV$  such that $\|\b_0\|_{\widetilde{E}_0(1,\infty)}<\infty$ and  $\|\b_1\|_{\widetilde{E}_1(0,1)}<\infty$ and denote
$$\rho(u) = u^{\theta_1- \theta_0} \frac{\a_0(u)\|\b_0\|_{\widetilde{E}_0(u,\infty)}}{\a_1(u)\|\b_1\|_{\widetilde{E}_1(0,u)} },\quad u>0.$$
\begin{itemize}
\item[a)] If $0<\theta<1$, then
$$
(\overline{X}^{\mathcal L}_{\theta_0, \b_0,E_0,\a_0,F_0}, \overline{X}^{\mathcal R}_{\theta_1,\b_1,E_1,\a_1,F_1})_{\theta,\b,E}= \overline{X}_{\tilde{\theta},\B_\theta,E}
$$
where $\tilde{\theta}=(1-\theta)\theta_0+\theta\theta_1$
and
$$B_\theta(u)=(\a_0(u)\|\b_0\|_{\widetilde{E}_0(u,\infty)})^{1-\theta}(\a_1(u)\|\b_1\|_{\widetilde{E}_1(0,u)})^{\theta}\b(\rho(u)),\quad u>0.$$ 

\item[b)] If $\theta=0$ and $\|\b\|_{\widetilde{E}(1,\infty)}<\infty$, then
$$(\overline{X}^{\mathcal L}_{\theta_0, \b_0,E_0,\a_0,F_0}, \overline{X}^{\mathcal R}_{\theta_1,\b_1,E_1,\a_1,F_1})_{0,\b,E}=\overline{X}^{\mathcal L}_{\theta_0,\B_0,E,\a_0,F_0}\cap\overline{X}^{\mathcal L,\mathcal L}_{\theta_0,\b\circ\rho,E,\b_0,E_0,a_0,F_0}$$
\vspace{1mm}
where $\B_0(u)=\|\b_0\|_{\widetilde{E}_0(u,\infty)}\b(\rho(u))$, $u>0$.

\item[c)]  If $\theta=1$ and  $\|\b\|_{\widetilde{E}(0,1)}<\infty$, then 
$$(\overline{X}^{\mathcal L}_{\theta_0, \b_0,E_0,\a_0,F_0}, \overline{X}^{\mathcal R}_{\theta_1,\b_1,E_1,\a_1,F_1})_{1,\b,E}=
\overline{X}^{\mathcal R}_{\theta_1,\B_1,E,\a_1,F_1}\cap\overline{X}^{\mathcal R,\mathcal R}_{\theta_1,\b\circ\rho,E,\b_1,E_1,a_1,F_1}$$
where $\B_1(u)=\|\b_1\|_{\widetilde{E}_1(0,u)}\b(\rho(u))$, $u>0$.
\end{itemize}
\end{thm}

\begin{proof}
As usual, we denote $Y_0=\overline{X}^{\mathcal L}_{\theta_0, \b_0,E_0,\a_0,F_0}$, $Y_1=\overline{X}^{\mathcal R}_{\theta_1,\b_1,E_1,\a_1,F_1}$ and $\overline{K}(u,f)=K(u,f;Y_0,Y_1)$, $u>0$. Lemma \ref{Le51} establishes the equivalence
$$\|f\|_{\overline{Y}_{\theta,\b,E}}=\|u^{-\theta}\b(u)\overline{K}(u,f)\|_{\widetilde{E}}\sim\|\rho(u)^{-\theta}\b(\rho(u))\overline{K}(\rho(u), f)\|_{\widetilde{E}}.$$
Then, in order to obtain a) it is sufficient to prove that
\begin{equation}\label{ec5RR}
\|\rho(u)^{-\theta}\b(\rho(u))\overline{K}(\rho(u), f)\|_{\widetilde{E}}\sim \big\|u^{- \tilde{\theta}} \B_\theta(u) K(u,f) \big \|_{\widetilde{E}}.
\end{equation}

Theorem \ref{5.17.2} and \eqref{e1} give
$$\overline{K}(\rho(u), f)\gtrsim
\|\b_0\|_{\widetilde{E}_0(u,\infty)}\| s^{- \theta_0} \a_0(s) K(s,f)\|_{ \widetilde{F}_0(0,u)}\gtrsim u^{- \theta_0} \a_0(u)\|\b_0\|_{\widetilde{E}_0(u,\infty)} K(u,f).$$
Hence, the inequality $``\gtrsim"$ of \eqref{ec5RR} follows using the relation
$$\rho(u)^{- \theta}  \b(\rho(u))=u^{\theta_0-\tilde{\theta}}\frac{\B_\theta(u)}{\a_0(u)\|\b_0\|_{\widetilde{E}_0(u,\infty)}},\quad u>0.$$
Next, we proceed with the reverse inequality $``\lesssim"$ in \eqref{ec5RR}. By Theorem \ref{5.17.2} and the triangular inequality, we have 
\begin{align*}
\big\|  \rho(u)^{-\theta} \b(\rho(u))  &\overline{K}(\rho(u),f)
\big\|_{\widetilde{E}} \\
& \lesssim 
 \Big \|\rho(u)^{-\theta} \b(\rho(u))  \big\|\b_0(t) 
     \|s^{-\theta_0}\a_0(s)K(s,f)\|_{\widetilde{F}_0(0,t)}\big\|_{\widetilde{E}_{0}(0,u)} \Big\|_{\widetilde{E}}\\
&+ \Big \|\rho(u)^{-\theta} \b(\rho(u)) \| \b_0\|_{\widetilde{E}_0(u,\infty)}
    \big\|t^{-\theta_0}\a_0(t)K(t,f)\big\|_{\widetilde{F}_0(0,u)}\Big\|_{\widetilde{E}}\\
&+  \Big\|\rho(u)^{1-\theta} \b(\rho(u)) \| \b_1\|_{\widetilde{E}_1(0,u)}
     \|t^{-\theta_1}\a_1(t)K(t,f)\|_{\widetilde{F}_1(u,\infty)}  \Big\|_{\widetilde{E}}\\
& + \Big\|\rho(u)^{1-\theta} \b(\rho(u)) \big\| \b_1(t) 
     \|s^{-\theta_1}\a_1(s)K(s,f)\|_{\widetilde{F}_1(t,\infty)}\big\|_{\widetilde{E}_{1}(u,\infty)}  \Big\|_{\widetilde{E}}\\
		&:=I_{14}+I_{15}+I_{16}+I_{17}. \end{align*}
Let us estimate the four expressions by $\|u^{- \tilde{\theta}} \B_\theta(u) K(u,f) \|_{\widetilde{E}}$. Lemma \ref{lema45}, with the monotone function $t\rightsquigarrow\|\cdot\|_{\widetilde{F}_0(0,t)}$, $\beta=\theta_0-\tilde{\theta}<0$ and Lemma \ref{lem1} (iv) yield that
\begin{align*}
I_{14}&\sim\Big \|  u^{\theta_0-\tilde{\theta}}\frac{B_\theta(u)}{\a_0(u)\|\b_0\|_{\widetilde{E}_0(u,\infty)}}\b_0(u)\|s^{-\theta_0}\a_0(s)K(s,f)\|_{\widetilde{F}_0(0,u)}\Big\|_{\widetilde{E}}\\
&\lesssim \Big \|  u^{\theta_0-\tilde{\theta}}\frac{B_\theta(u)}{\a_0(u)}\|s^{-\theta_0}\a_0(s)K(s,f)\|_{\widetilde{F}_0(0,u)}\Big\|_{\widetilde{E}}.
\end{align*} 
Lemma \ref{lema45} again, with the $K$-functional, implies that $I_{14}\lesssim\|u^{- \tilde{\theta}} \B_\theta(u) K(u,f) \|_{\widetilde{E}}.$ Similarly, using Lemma \ref{lema45} only once, we establish that $I_{15}\sim\|u^{- \tilde{\theta}} \B_\theta(u) K(u,f) \|_{\widetilde{E}}$. Having in mind that
\begin{equation*}\label{904}
\rho(u)^{1- \theta}  \b(\rho(u))=u^{\theta_1-\tilde{\theta}}\frac{\B_\theta(u)}{\a_1(u)\|\b_1\|_{\widetilde{E}_1(0,u)}},\quad u>0
\end{equation*}
and arguing as we did with $I_2$ in Theorem \ref{thmRR}, we have that 
$I_{16}\sim\|u^{\tilde{\theta}}\B_\theta(u)K(u,f)\|_{\widetilde{E}}.$
The estimate of the term $I_{17}$ can be carried out, arguing as in $I_{14}$, using Lemma \ref{lema45} twice, once with the monotone function $t\rightsquigarrow\|\cdot\|_{\widetilde{F}_1(t,\infty)}$, the second one with the $K$-functional, $\gamma=\theta_1-\tilde{\theta}>0$ and Lemma \ref{lem1} (iv). 
In fact,
\begin{align*}
I_{17}&=\Big\| u^{\theta_1-\tilde{\theta}}\frac{\B_\theta(u)}{\a_1(u)\|\b_1\|_{\widetilde{E}_1(0,u)}}\big\| \b_1(t) 
     \|s^{-\theta_1}\a_1(s)K(s,f)\|_{\widetilde{F}_1(t,\infty)}\big\|_{\widetilde{E}_{1}(u,\infty)}  \Big\|_{\widetilde{E}}\\
&\sim\Big\| u^{\theta_1-\tilde{\theta}}\frac{\B_\theta(u)}{\a_1(u)\|\b_1\|_{\widetilde{E}_1(0,u)}} \b_1(t) 
     \|s^{-\theta_1}\a_1(s)K(s,f)\|_{\widetilde{F}_1(t,\infty)} \Big\|_{\widetilde{E}}\\	
&\lesssim\Big\| u^{\theta_1-\tilde{\theta}}\frac{\B_\theta(u)}{\a_1(u)} 
     \|s^{-\theta_1}\a_1(s)K(s,f)\|_{\widetilde{F}_1(t,\infty)} \Big\|_{\widetilde{E}}\\	
		&\sim\Big\|u^{\tilde{\theta}}\B_\theta(u)K(u,f)\Big\|_{\widetilde{E}}.	
\end{align*}		
This concludes the proof of a). 

To prove b) it suffices to show that 
\begin{equation}\label{e99}
\|\b(\rho(u))\overline{K}(\rho(u), f)\|_{\widetilde{E}}\sim \max\big\{\|f\|_{\overline{X}^{\mathcal L}_{\theta_0,\B_0,E,\a_0,F_0}},\|f\|_{\overline{X}^{\mathcal L,\mathcal L}_{\theta_0,\b\circ\rho,E,\b_0,E_0,a_0,F_0}}\big\}.
\end{equation}
Theorem \ref{5.17.2} guarantees that 
$$\overline{K}(\rho(u),f)\gtrsim\|\b_0\|_{\widetilde{E}_0(u,\infty) }\|t^{- \theta_0} \a_0(t) K(t,f)\|_{ \widetilde{F}_0(0,u)}$$
and 
$$\overline{K}(\rho(u),f)\gtrsim \Big\| \b_0(t) \|s^{- \theta_0} \a_0(s) K(s,f)\|_{ \widetilde{F}_0(0,t)}\Big\|_{ \widetilde{E}_0(0,u)}.$$
Hence the inequality ``$\gtrsim$'' of \eqref{e99} follows. On the other hand, Theorem \ref{5.17.2} and triangular inequality give that
\begin{align*}
\big\| \b(\rho(u))  \overline{K}(\rho(u),f)
\big\|_{\widetilde{E}} 
& \lesssim 
 \Big \| \b(\rho(u))  \big\|\b_0(t) 
     \|s^{-\theta_0}\a_0(s)K(s,f)\|_{\widetilde{F}_0(0,t)}\big\|_{\widetilde{E}_{0}(0,u)} \Big\|_{\widetilde{E}}\\
&+ \Big \| \b(\rho(u)) \| \b_0\|_{\widetilde{E}_0(u,\infty)}
    \big\|t^{-\theta_0}\a_0(t)K(t,f)\big\|_{\widetilde{F}_0(0,u)}\Big\|_{\widetilde{E}}\\
&+  \Big\|\rho(u) \b(\rho(u)) \| \b_1\|_{\widetilde{E}_1(0,u)}
     \|t^{-\theta_1}\a_1(t)K(t,f)\|_{\widetilde{F}_1(u,\infty)}  \Big\|_{\widetilde{E}}\\
& + \Big\|\rho(u)\b(\rho(u)) \big\| \b_1(t) 
     \|s^{-\theta_1}\a_1(s)K(s,f)\|_{\widetilde{F}_1(t,\infty)}\big\|_{\widetilde{E}_{1}(u,\infty)}  \Big\|_{\widetilde{E}}\\
		&:=I_{18}+I_{19}+I_{20}+I_{21}. \end{align*}
It is clear that $I_{18}=\|f\|_{\overline{X}^{\mathcal L,\mathcal L}_{\theta_0,\b\circ\rho,E,\b_0,E_0,a_0,F_0}}$ and $I_{19}=\|f\|_{\overline{X}^{\mathcal L}_{\theta_0,\B_0,E,\a_0,F_0}}$. Besides, Lemma \ref{lema45}, with the $K$-functional and $\gamma=\theta_1-\theta_0>0$, and inequality \eqref{e1} yield 
\begin{align*}
I_{20}&=\Big\|u^{\theta_1-\theta_0}\frac{\a_0(u) \| \b_0\|_{\widetilde{E}_0(u,\infty)}}{\a_1(u)}\b(\rho(u)) 
     \|t^{-\theta_1}\a_1(t)K(t,f)\|_{\widetilde{F}_1(u,\infty)}  \Big\|_{\widetilde{E}}\\
		&\sim \Big\|u^{-\theta_0}\a_0(u) \| \b_0\|_{\widetilde{E}_0(u,\infty)}\b(\rho(u)) K(u,f) \Big\|_{\widetilde{E}}\\
	&\lesssim 	\Big\|\| \b_0\|_{\widetilde{E}_0(u,\infty)}\b(\rho(u))\|t^{- \theta_0} \a_0(t) K(t,f)\|_{ \widetilde{F}_0(0,u)}\Big\|_{\widetilde{E}}=\|f\|_{\overline{X}^{\mathcal L}_{\theta_0,\B_0,E,\a_0,F_0}}.
\end{align*}
The estimate of the last term $I_{21}$ requires, as with $I_{17}$, the use of Lemma \ref{lema45} twice, with $\gamma=\theta_1-\theta_0>0$, and Lemma \ref{lem1} (iv). Indeed,
\begin{align*}
I_{21}&=\Big\|u^{\theta_1-\theta_0}\frac{\a_0(u) \| \b_0\|_{\widetilde{E}_0(u,\infty)}}{\a_1(u)\|\b_1\|_{\widetilde{E}_1(0,u)}}\b(\rho(u)) \big\| \b_1(t) 
     \|s^{-\theta_1}\a_1(s)K(s,f)\|_{\widetilde{F}_1(t,\infty)}\big\|_{\widetilde{E}_{1}(u,\infty)}  \Big\|_{\widetilde{E}}\\
		&\sim\Big\|u^{\theta_1-\theta_0}\frac{\a_0(u) \| \b_0\|_{\widetilde{E}_0(u,\infty)}}{\a_1(u)\|\b_1\|_{\widetilde{E}_1(0,u)}}\b(\rho(u))  \b_1(u) 
     \|s^{-\theta_1}\a_1(s)K(s,f)\|_{\widetilde{F}_1(u,\infty)}  \Big\|_{\widetilde{E}}\\
		&\lesssim \Big\|u^{\theta_1-\theta_0}\frac{\a_0(u) \| \b_0\|_{\widetilde{E}_0(u,\infty)}}{\a_1(u)}\b(\rho(u)) 
     \|s^{-\theta_1}\a_1(s)K(s,f)\|_{\widetilde{F}_1(u,\infty)}  \Big\|_{\widetilde{E}}\\
		&\sim \Big\|u^{-\theta_0}\a_0(u) \| \b_0\|_{\widetilde{E}_0(u,\infty)}\b(\rho(u))  
     K(u,f)  \Big\|_{\widetilde{E}}.
		\end{align*}
Finally, using \eqref{e1} we have that $I_{21}\lesssim \|f\|_{\overline{X}^{\mathcal L}_{\theta_0,\B_0,E,\a_0,F_0}}$ . This establishes b).

The proof of c) can be done using similar arguments to Theorem \ref{thmRR} c).
\end{proof}

\section{Applications}\label{applications}

For simplicity, we apply our results to ordered (quasi)-Banach
couples $X = (X_0, X_1)$, in the sense that  $X_1\hookrightarrow X_0$. The most classical example of an ordered couple is $(L_1(\Omega,\mu),L_\infty(\Omega,\mu))$ when $\Omega$ is a finite measure space. 

\subsection{Ordered couples}\label{aordered}\hspace{2mm}

\vspace{2mm}
We briefly review how our definitions adapt to this simpler setting of ordered couples.
When $X_1\hookrightarrow X_0$, the real interpolation  $\overline{X}_{\theta,\b,E}$ can be equivalently defined as
$$\overline{X}_{\theta,\b,E}=\big\{f\in X_0: \|t^{-\theta}\b(t)K(t,f)\|_{\widetilde{E}(0,1)}<\infty\big\}$$
where $0\leq \theta\leq 1$, $E$ is an r.i. space on $(0,1)$ and $\b\in SV(0,1)$.  The replacement of $(0,\infty)$ by $(0,1)$ in Definition \ref{def1} yields the definition of the class $SV(0, 1)$.

Similarly, given a real parameter $0 \leq \theta \leq 1$,  $a, \b, c\in SV(0,1)$ and r.i. spaces $E, F, G$ on $(0,1)$, the spaces  
$\overline{X}_{\theta,\b,E,a,F}^{\mathcal L}$, $\overline{X}_{\theta,c,E,\b,F,a,G}^{\mathcal L,\mathcal L }$, $\overline{X}_{\theta,c,E,\b,F,a,G}^{\mathcal R,\mathcal L }$ are defined just as in Definitions \ref{defLR} and \ref{defLRR}; the only change being that $\widetilde{E}(0,\infty)$ must be replaced by $\widetilde{E}(0, 1)$.
Likewise, the spaces  $\overline{X}_{\theta,\b,E,a,F}^{\mathcal R}$, 
$\overline{X}_{\theta,c,E,\b,F,a,G}^{\mathcal R,\mathcal R}$ and $\overline{X}_{\theta,c,E,\b,F,a,G}^{\mathcal L,\mathcal R}$ are defined as
$$\overline{X}_{\theta,\b,E,a,F}^{\mathcal R}=\Big\{f\in X_0:\big \|  \b(t) \|   s^{-\theta} a(s) K(s,f) \|_{\widetilde{F}(t,1)}      \big   \|_{\widetilde{E}(0,1)}<\infty\Big\},$$
$$\overline{X}_{\theta,c,E,\b,F,a,G}^{\mathcal R,\mathcal R }=\bigg\{f\in X_0: 
\Big\|c(u)\big\|\b(t) \|s^{-\theta} \a(s) K(s,f) \|_{\widetilde{G}(t,1)}\big\|_{\widetilde{F}(u,1)}\Big\|_{\widetilde{E}(0,1)} < \infty\bigg\}$$
and
$$\overline{X}_{\theta,c,E,\b,F,a,G}^{\mathcal L,\mathcal R }=\bigg\{f\in X_0: 
\Big\|c(u)\big\|\b(t) \|s^{-\theta} \a(s) K(s,f) \|_{\widetilde{G}(u,t)}\big\|_{\widetilde{F}(u,1)}\Big\|_{\widetilde{E}(0,1)} < \infty\bigg\}$$

Of course, all the results in the paper remain true if we work with an ordered couple  and use as parameters slowly varying functions on $(0,1)$ and r.i. spaces on $(0,1)$. In these cases all assumptions concerning the interval $(1,\infty)$ must be omitted.

Moreover, if the couple is ordered, then the ${\mathcal R}$ and ${\mathcal L}$-scales  are also ordered.
\begin{lem}\label{lemainclusion}
Let $\overline{X}$ be an ordered (quasi-) Banach couple,  $E_0,\ E_1,\ F_0,\ F_1$ r.i. spaces on $(0,1)$ and $a_0, a_1, \b_0, \b_1\in SV(0,1)$. If $0<\theta_0<\theta_1<1$, then
$$\overline{X}^{\mathcal R} _{\theta_1,\b_1,E_1,a_1,F_1}\hookrightarrow \overline{X}^{\mathcal R}_{\theta_0,\b_0,E_0,a_0,F_0}, \quad \quad \quad  \overline{X}^{\mathcal L} _{\theta_1,\b_1,E_1,a_1,F_1}\hookrightarrow \overline{X}^{\mathcal L}_{\theta_0,\b_0,E_0,a_0,F_0},$$
$$\overline{X}^{\mathcal R} _{\theta_1,\b_1,E_1,a_1,F_1}\hookrightarrow \overline{X}^{\mathcal L}_{\theta_0,\b_0,E_0,a_0,F_0},\quad  \quad \quad\overline{X}^{\mathcal L} _{\theta_1,\b_1,E_1,a_1,F_1}\hookrightarrow \overline{X}^{\mathcal R}_{\theta_0,\b_0,E_0,a_0,F_0},$$
assuming, if it is necessary, the condition $\|\b_0\|_{\widetilde{E}_0(0,1)}<\infty$ such that the right hand side of each inclusion is not the trivial space.
\end{lem}

\begin{proof}
Let $f\in \overline{X}^{\mathcal R}_{\theta_1,\b_1,E_1,a_1,F_1}$ and assume that $\|\b_0\|_{\widetilde{E}_0(0,1)}<\infty$ and $\|\b_1\|_{\widetilde{E}_1(0,1)}<\infty$, otherwise the first inclusion is trivial. Using \eqref{eK2}, that the function $t\rightsquigarrow\|\cdot\|_{\widetilde{F}_0(t,1)}$ is non-increasing and Lemma \ref{lem1} (i) with $\alpha=\theta_1-\theta_0>0$, we obtain
\begin{align*}
\|f\|_{{\mathcal R};\theta_0,\b_0,E_0,a_0,F_0}
&\lesssim \|f\|_{{\mathcal R};\theta_1,\b_1,E_1,a_1,F_1}\bigg \|\b_0(t)\Big\|s^{\theta_1-\theta_0}\frac{a_0(s)}{a_1(s)\|\b_1\|_{\widetilde{E}_1(0,s)}} \Big\|_{\widetilde{F}_0(t,1)}      \bigg   \|_{\widetilde{E}_0(0,1)} \\
&\leq \|f\|_{{\mathcal R};\theta_1,\b_1,E_1,a_1,F_1}\|\b_0\|_{\widetilde{E}_0(0,1)}\Big\|s^{\theta_1-\theta_0}\frac{a_0(s)}{a_1(s)\|\b_1\|_{\widetilde{E}_1(0,s)}} \Big\|_{\widetilde{F}_0(0,1)}   \\
&\sim \|f\|_{{\mathcal R};\theta_1,\b_1,E_1,a_1,F_1} 
\end{align*}
which gives that $f\in \overline{X}^{\mathcal R}_{\theta_0,\b_0,E_0,a_0,F_0}$. The proof of the other three inclusions can be carried out using similar arguments, using \eqref{eK3} in place of \eqref{eK2} if it is necessary.
\end{proof}
\subsection{ Interpolation between grand and small Lebesgue spaces}\hspace{2mm}

\vspace{2mm}
Let $(\Omega,\mu)$ denote totally $\sigma$-finite measure space and $\mathcal{M}(\Omega,\mu)$ the set of all $\mu$-measurable functions on $(\Omega,\mu)$.

\begin{defn}\label{ultrasymmetric}
Let $1<p\leq\infty$, $\b\in SV$ and $E$ an r.i. space. The \textit{Lorentz-Karamata type} space $L_{p,\b,E}$ is defined as the set of all $f\in\mathcal{M}(\Omega,\mu)$ such that
$$\|f\|_{L_{p,\b,E}}=\|t^{1/p}\b(t)f^*(t)\|_{\widetilde{E}}<\infty.$$
\end{defn}
When $p=\infty$, $L_{p,\b,E}$ is not the trivial space if, and only if, $\|\b\|_{\widetilde{E}(0,1)}<\infty$.

These spaces are particular examples of the ultrasymmetric spaces studied by E. Pustylnik \cite{Pust-ultrasymm}. When $E=L_q$, $1\leq q\leq \infty$, the space coincides with the classical \textit{Lorentz-Karamata} space $L_{p,q;\b}$ (see \cite{GOT,neves}). 

Peetre's well-known formula \cite{Bennett-Sharpley,Peetre}
$$K(t,f;L_1,L_\infty)=\int_0^tf^*(s)\, ds= tf^{**}(t),\quad t>0,$$
and the equivalence $\|t^{1/p}\b(t)f^{**}(t)\|_{\widetilde{E}}\sim \|t^{1/p}\b(t)f^*(t)\|_{\widetilde{E}}$ , $1<p\leq \infty$   (see, e.g. \cite[Lemma 2.16]{CwP}), yield that
\begin{equation}\label{eU}
(L_1,L_\infty)_{1-\frac1p,\b,E}=L_{p,\b,E}.
\end{equation}
Analogously it can be proved that
$$(L_1,L_\infty)^{\mathcal R}_{1-\frac1p,\b,E,a,F}=L^{\mathcal R}_{p,\b,E,a,F}\mand (L_1,L_\infty)^{\mathcal L}_{1-\frac1p,\b,E,a,F}=L^{\mathcal L}_{p,\b,E,a,F}$$
where
$$L^{\mathcal R}_{p,\b,E,a,F}:=\{f\in\mathcal{M}(\Omega,\mu):
\Big\|\b(t) \|s^{1/p} a(s) f^*(s)\|_{\widetilde{F}(t,\infty)}\Big\|_{\widetilde{E}} < \infty\}$$
and 
$$L^{\mathcal L}_{p,\b,E,a,F}:=\{f\in\mathcal{M}(\Omega,\mu):
\Big\|\b(t) \|s^{1/p} a(s) f^*(s)\|_{\widetilde{F}(0,t)}\Big\|_{\widetilde{E}} < \infty\}$$
for $E$, $F$ r.i. spaces, $a$, $\b\in SV$ and $1<p\leq \infty$
(see \cite{GOT} for the case $E=L_r$, $F=L_q$, $0<q,r\leq \infty$). Also, in a similar way
$$(L_1,L_\infty)^{\mathcal R,\mathcal R}_{1-\frac1p,c,E,\b,F,a,G}=L^{\mathcal R,\mathcal R}_{p,c, E,\b,F,a,G},\qquad (L_1,L_\infty)^{\mathcal L,\mathcal L}_{1-\frac1p,c,E,\b,F,a,G}=L^{\mathcal L,\mathcal L}_{p,c,E,\b,F,a,G},$$
$$(L_1,L_\infty)^{\mathcal R,\mathcal L}_{1-\frac1p,c,E,\b,F,a,G}=L^{\mathcal R,\mathcal L}_{p,c,E,\b,F,a,G}\mand (L_1,L_\infty)^{\mathcal L,\mathcal R}_{1-\frac1p,c,E,\b,F,a,G}=L^{\mathcal L,\mathcal R}_{p,c,E,\b,F,a,G}$$
where the spaces $L^{\mathcal R,\mathcal R}_{p,c, E,\b,F,a,G}$, $L^{\mathcal L,\mathcal L}_{p,c,E,\b,F,a,G}$, $L^{\mathcal R,\mathcal L}_{p,c,E,\b,F,a,G}$ and $L^{\mathcal L,\mathcal R}_{p,c,E,\b,F,a,G}$ are defined as the set of all $f\in\mathcal{M}(\Omega,\mu)$ for which \eqref{dRR}-\eqref{dLR} are satisfied changing  $s^{-\theta}a(s)K(s,f)$ by $s^{1/p}a(s)f^*(s)$, $s>0$, respectively.

Following the paper by Fiorenza and Karadzhov \cite{FK} we give the next definition:
\begin{defn}\label{def_gLp}
Let $(\Omega,\mu)$ be a finite measure space such that $\mu(\Omega)=1$, let $1<p<\infty$ and $\alpha>0$. The grand Lebesgue space $L^{p),\alpha}$ is  the set of all $f\in{\mathcal M}(\Omega,\mu)$ such that
\begin{equation*}\label{lpa1}
\|f\|_{p),\alpha}=\Big\|\ell^{-\frac{\alpha}{p}}(t)\|s^{1/p}f^{*}(s)\|_{\widetilde{L}_p(t,1)}\Big\|_{L_\infty(0,1)}<\infty.
\end{equation*}
The small Lebesgue space $L^{(p,\alpha}(\Omega)$ is  the set of all $f\in{\mathcal M}(\Omega,\mu)$ such that
\begin{equation*}\label{lpa2}
\|f\|_{(p,\alpha}=\Big\|\ell^{\frac{\alpha}{p'}-1}(t)\|s^{1/p}f^{*}(s)\|_{\widetilde{L}_p(0,t)}\Big\|_{\widetilde{L}_1(0,1)}<\infty
\end{equation*}
where $\frac1p+\frac1{p'}=1$.
\end{defn}

The classical grand Lebesgue space $L^{p)}(\Omega):=L^{p),1}(\Omega)$ was introduced by Iwaniec and Sbordone in \cite{IS} while the classical small Lebesgue space $L^{(p}(\Omega):=L^{(p,1}(\Omega)$ was characterized by Fiorenza  in \cite{F1} as its associate; that is $ (L^{(p'})' = L^{p)}.$ For more information about this spaces and their generalizations see the recent paper \cite{FFG}.

As observed in \cite{FK,op1}, these spaces can be characterized as $\mathcal{R}$ and  $\mathcal{L}$-spaces in the following way
$$L^{p),\alpha}=L^{\mathcal R}_{p,\ell^{-\alpha/p}(u),L_\infty,1,L_p}\mand
L^{(p,\alpha}=L^{\mathcal L}_{p,\ell^{\alpha/p'-1}(u),L_1,1,L_p}.$$

Thus, we can apply the results from \S\ref{reiterationtheorems} to identify the interpolation of grand and small Lebesgue spaces. Moreover, the following technical lemma is needed.
\begin{lem}\cite[Lemma 6.1]{EOP}
If $\sigma+\frac{1}{q}<0$  with $1\leq q<\infty$ or $q=\infty$ and $\sigma\leq 0$, then 
\begin{equation}\label{einfty1}
\|\ell^{\sigma}(t)\|_{\widetilde{L}_{q}(0,u)}\sim \ell^{\sigma+\frac1{q}}(u), \quad u\in(0,1).
\end{equation}
If $\sigma+\frac{1}{q}>0$  with $1\leq q<\infty$, or $q=\infty$ and $\sigma\geq 0$, then
 \begin{equation}\label{einfty2}
\|\ell^{\sigma}(t)\|_{\widetilde{L}_{q}(u,1)}\sim \ell^{\sigma+\frac1{q}}(u),\quad u\in(0,1/2).\end{equation}
\end{lem}

\vspace{2mm}
Notice also that if $\b(t)\sim\a(t)$ for all $t\in(0,1/2)$, then the monotonicity properties of the $K$-functional and the properties of the slowly varying functions imply that
$$\overline{X}_{\theta,\b,E}=\overline{X}_{\theta,\a,E}.$$
Thus, for any $0<\theta<1$, $a\in SV(0,1)$ and r.i. space $E$ on $(0,1)$,
\begin{equation*}
\overline{X}_{\theta,\|\ell^{\sigma}(t)\|_{\widetilde{L}_{q}(u,1)},E}=\overline{X}_{\theta, \ell^{\sigma+\frac1{q}}(u),E}.
\end{equation*}
A similar identity holds for $\mathcal{R}$, $\mathcal{L}$-spaces and the extreme constructions.

Using Theorem \ref{thmRR} and \eqref{einfty1} with $q=\infty$, we  can state the interpolation formulae for a couple formed by two grand Lebesgue spaces. Similar results in the no-limiting cases $0<\theta<1$ appear in \cite[Theorem 8.3]{AFH}, \cite[Corollary 51]{Do2020} and \cite[Theorem 1.2]{FFGKR}. The result also completes \cite[Corollary 5.7]{FMS-RL1}.

\begin{cor}\label{cor55}
Let $E$ be an r.i. space on $(0,1)$, $\b\in SV(0,1)$, $1< p_0<p_1<\infty$ and $\alpha, \beta>0$. Consider the function
$\rho(u)=u^{\frac1{p_0}-\frac1{p_1}}\ell^{\frac{\beta}{p_1}-\frac{\alpha}{p_0}}(u)$, $u\in (0,1)$.

\begin{itemize}
\vspace{1mm}
\item[a)] If $0<\theta< 1$, then
\begin{equation*}
\big(L^{p_0),\alpha},L^{p_1),\beta}\big)_{\theta,\b,E}=L_{p,\B_\theta,E},
\end{equation*}
where $\frac1{p}=\frac{1-\theta}{p_0}+\frac{\theta}{p_1}$ and $\B_\theta(u)=\ell^{-\big[\frac{\alpha(1-\theta)}{p_0}+\frac{\beta\theta}{p_1}\big]}(u)\b(\rho(u))$, $u \in(0,1)$.

\vspace{1mm}
\item[b)] If $\theta=0$, then
\begin{equation*}
\big(L^{p_0),\alpha},L^{p_1),\beta}\big)_{0,\b,E}=L_{p_0,\b\circ\rho,E,\ell^{-\alpha/p_0}(u),L_\infty,1,L_{p_0}}^{\mathcal R,\mathcal L}.
\end{equation*}

\vspace{1mm}
\item[c)] If $\theta=1$ and $\|\b\|_{\widetilde{E}(0,1)}<\infty$, then
\begin{equation}\label{einter1}
\big(L^{p_0),\alpha},L^{p_1),\beta}\big)_{1,\b,E}=L^{\mathcal R}_{p_1,\B_1,E,1,L_{p_1}}\cap L_{p_1,\b\circ\rho,E,\ell^{-\beta/p_1}(u),L_\infty,1,L_{p_1}}^{\mathcal R,\mathcal R}
\end{equation}
where $\B_1(u)=\ell^{\frac{-\beta}{p_1}}(u)\b(\rho(u))$,  $u \in(0,1)$.
\end{itemize}
\end{cor}

\begin{rem}
Comparing \eqref{einter1} with Corollary 5.7 c) from \cite{FMS-RL1}, we observe that
$$L^{\mathcal R}_{p_1,\B_1,E,1,L_{p_1}}\cap(L^{p_0),\alpha},L_{p_1})^{\mathcal R}_{1,\b\circ\rho^{\#},E,\ell^{-\beta/p_1}(u),L_\infty} =L^{\mathcal R}_{p_1,\B_1,E,1,L_{p_1}}\cap L_{p_1,\b\circ\rho,E,\ell^{-\beta/p_1}(u),L_\infty,1,L_{p_1}}^{\mathcal R,\mathcal R}$$
where $\rho^{\#}(u)=u\ell^{\frac{\beta}{p_1}}(u)$, $u\in(0,1)$.
\end{rem}

Theorem \ref{thmLL} and  \eqref{einfty2} with $q=1$, enables us  to state the  following interpolation formulae for a couple formed by two small Lebesgue spaces. The result completes \cite[Corollary 49]{Do2020}, \cite[Corollary 5.12]{FMS-RL1} and \cite[Theorem 3.4]{FFGKR}. 

\begin{cor}\label{cor57}
Let $E$ be an r.i. space on $(0,1)$, $\b\in SV(0,1)$, $1< p_0<p_1<\infty$ and $\alpha, \beta>0$. Consider the function
$\rho(u)=u^{\frac1{p_0}-\frac1{p_1}}\ell^{\frac{\alpha}{p'_0}-\frac{\beta}{p'_1}}(u)$, $u\in (0,1)$.

\begin{itemize}
\vspace{1mm}
\item[a)] If $0<\theta< 1$, then
\begin{equation*}
\big(L^{(p_0,\alpha},L^{(p_1,\beta}\big)_{\theta,\b,E}=L_{p,\B_\theta,E},
\end{equation*}
where $\frac1{p}=\frac{1-\theta}{p_0}+\frac{\theta}{p_1}$ and $\B_\theta(u)=\ell^{\frac{\alpha(1-\theta)}{p'_0}+\frac{\beta\theta}{p'_1}}(u)\b(\rho(u))$, $u \in(0,1)$.

\vspace{1mm}
\item[b)] If $\theta=0$, then
\begin{equation}\label{einter2}
\big(L^{(p_0,\alpha},L^{(p_1,\beta}\big)_{0,\b,E}=L^{\mathcal L}_{p_0,B_0,E,1,L_{p_0}}\cap L^{\mathcal L,\mathcal L}_{p_0,\b\circ\rho,E,\ell^{\alpha/p'_0-1}(u),L_1,1,L_{p_0}}
\end{equation}
where $B_0(u)=\ell^{\frac{\alpha}{p'_0}}(u)\b(\rho(u))$, $u\in (0,1)$.
\vspace{1mm}
\item[c)] If $\theta=1$ and $\|\b\|_{\widetilde{E}(0,1)}<\infty$, then
$$
\big(L^{(p_0,\alpha},L^{(p_1,\beta}\big)_{1,\b,E}=L^{\mathcal L,\mathcal R}_{p_1,\b\circ\rho,E,\ell^{\beta/p'_1-1}(u),L_1,1,L_{p_1}}.$$
\end{itemize}
\end{cor}

\begin{rem} 
Comparing \eqref{einter2} with Corollary 5.12 b) from \cite{FMS-RL1}, we observe that
$$L^{\mathcal L}_{p_0,B_0,E,1,L_{p_0}}\cap\big(L_{p_0},L^{(p_1,\beta}\big)^{\mathcal L}_{0,\b\circ\rho^{\#},E,\ell^{\alpha/p'_0}(u),L_1}=
L^{\mathcal L}_{p_0,B_0,E,1,L_{p_0}}\cap L^{\mathcal L,\mathcal L}_{p_0,\b\circ\rho,E,\ell^{\alpha/p'_0-1}(u),L_1,1,L_{p_0}}
$$
where $\rho^\#(u)=u\ell^{\alpha/p'_0}(u)$, $u\in (0,1)$.
\end{rem}

Finally, using Theorems \ref{thmRL}, \ref{thmLR} and estimates \eqref{einfty1}, \eqref{einfty2}, we identify the interpolation space between grand and small Lebesgue spaces. Corollary \ref{cor1} completes \cite[Corollary 52]{Do2020} with the limiting cases $\theta=0,1$, while  Corollary \ref{cor2} completes  \cite[Corollary 50]{Do2020}, \cite[Theorem 5.7]{FMS-RL2} and \cite[Theorem 5.7]{FFGKR}.

\begin{cor}\label{cor1}
Let $E$ be an r.i. space on $(0,1)$, $\b\in SV(0,1)$, $1< p_0<p_1<\infty$ and $\alpha, \beta>0$. Consider the function
$\rho(u)=u^{\frac1{p_0}-\frac1{p_1}}\ell^{-[\frac{\alpha}{p_0}+\frac{\beta}{p'_1}]}(u)$, $u\in (0,1)$.

\begin{itemize}

\vspace{1mm}
\item[a)] If $0<\theta< 1$, then
\begin{equation*}
\big(L^{p_0),\alpha},L^{(p_1,\beta}\big)_{\theta,\b,E}=L_{p,\B_\theta,E},
\end{equation*}
where $\frac1{p}=\frac{1-\theta}{p_0}+\frac{\theta}{p_1}$ and $\B_\theta(u)=\ell^{-\frac{\alpha(1-\theta)}{p_0}+\frac{\beta\theta}{p'_1}}(u)\b(\rho(u))$, $u \in(0,1)$.

\vspace{1mm}
\item[b)] If $\theta=0$, then
\begin{equation*}
\big(L^{p_0),\alpha},L^{(p_1,\beta}\big)_{0,\b,E}=L^{\mathcal R,\mathcal L}_{p_0,\b\circ\rho,E,\ell^{-\alpha/p_0}(u),L_\infty,1,L_{p_0}}.
\end{equation*}
\vspace{1mm}
\item[c)] If $\theta=1$ and $\|\b\|_{\widetilde{E}(0,1)}<\infty$, then
\begin{equation*}
\big(L^{p_0),\alpha},L^{(p_1,\beta}\big)_{1,\b,E}=L^{\mathcal L,\mathcal R}_{p_1,\b\circ\rho,E,\ell^{\beta/p'_1-1}(u),L_1,1,L_{p_1}}.
\end{equation*}
\end{itemize}
\end{cor}

\begin{cor}\label{cor2}
Let $E$ be an r.i. space on $(0,1)$, $\b\in SV(0,1)$, $1< p_0<p_1<\infty$ and $\alpha, \beta>0$. Consider the function
$\rho(u)=u^{\frac1{p_0}-\frac1{p_1}}\ell^{\frac{\alpha}{p'_0}+\frac{\beta}{p_1}}(u)$, $u\in (0,1)$.
\begin{itemize}

\vspace{1mm}
\item[a)] If $0<\theta< 1$, then
\begin{equation*}
\big(L^{(p_0,\alpha},L^{p_1),\beta}\big)_{\theta,\b,E}=L_{p,\B_\theta,E},
\end{equation*}
where $\frac1{p}=\frac{1-\theta}{p_0}+\frac{\theta}{p_1}$ and $\B_\theta(u)=\ell^{\frac{\alpha(1-\theta)}{p'_0}-\frac{\beta\theta}{p_1}}(u)\b(\rho(u))$, $u \in(0,1)$.

\vspace{1mm}
\item[b)] If $\theta=0$, then
\begin{equation}\label{einter3}
\big(L^{(p_0,\alpha},L^{p_1),\beta}\big)_{0,\b,E}=L^{\mathcal L}_{p_0,B_0,E,1,L_{p_0}}\cap L^{\mathcal L,\mathcal L}_{p_0,\b\circ\rho,E,\ell^{\alpha/p'_0-1}(u),L_1,1,L_{p_0}}
\end{equation}
where $B_0(u)=\ell^{\frac{\alpha}{p'_0}}(u)\b(\rho(u))$, $u\in (0,1)$.

\vspace{1mm}
\item[c)] If $\theta=1$ and $\|\b\|_{\widetilde{E}(0,1)}<\infty$, then
\begin{equation*}
\big(L^{(p_0,\alpha},L^{p_1),\beta}\big)_{1,\b,E}
=L^{\mathcal R}_{p_1,\B_1,E,1,L_{p_1}}\cap L^{\mathcal R,\mathcal R}_{p_1,\b\circ\rho,E,\ell^{-\beta/p_1}(u),L_\infty,1,L_{p_1}}
\end{equation*}
where $\B_1(u)=\ell^{-\frac{\beta}{p_1}}(u)\b(\rho(u))$, $u \in(0,1)$.

\end{itemize}
\end{cor}

\begin{rem}
Comparing \eqref{einter3} with Theorem 5.7 b) from \cite{FMS-RL2} one can deduce the identity
$$L^{\mathcal L}_{p_0,B_0,E,1,L_{p_0}}\cap(L_{p_0},L^{p_1),\beta})^{\mathcal L}_{0,b\circ\rho^{\#},E,\ell^{\alpha/p'_0-1}(u),L_1}=L^{\mathcal L}_{p_0,B_0,E,1,L_{p_0}}\cap L^{\mathcal L,\mathcal L}_{p_0,\b\circ\rho,E,\ell^{\alpha/p'_0}(u),L_1,1,L_{p_0}}$$
where $\rho^{\#}(u)=u\ell^{\frac{\alpha}{p'_0}}(u)$, $u\in(0,1)$. 
\end{rem}

\subsection{Generalized Gamma spaces.}\hspace{1mm}

\vspace{2mm}
Our goal in this subsection is to obtain interpolation formulae for couples formed by two Generalized Gamma spaces with double weight; see \cite{FFGKR}. 

\begin{defn}
Let $1 \leq p < \infty$, $1 \leq q \leq \infty$ and $w_{1}$, $w_{2}$ two weights on $(0,1)$ satisfying the following conditions:

\vspace{2mm}
\begin{enumerate}
\item[(c1)] There exist $K_{21}>0$ such that $w_2(2t)\leq K_{12}w_2(t)$,  for all $t\in(0,1/2)$. The space $L^p(0,1;w_2)$ is continuously embedded in $L^1(0,1)$.
\item[(c2)] The function $\int_0^t w_2(s)ds$ belongs to $L^{\frac{q}{p}}(0,1;w_1).$
\end{enumerate}
\vspace{2mm}The \textit{Generalized Gamma} space with double weights $G\Gamma(p,q,w_1,w_2)$ is the set of all measurable functions $f$ on $(0,1)$ such that 
$$\|f\|_{G\Gamma}=\bigg(\int_0^1w_1(t)\Big(\int_0^tw_2(s)(f^*(s))^p\ ds\Big)^{\frac{q}{p}}dt\bigg)^{\frac{1}{q}}<\infty.$$
\end{defn}

These spaces are a generalization of the spaces $G\Gamma(p,q,w_1):= G\Gamma(p,q,w_1,1)$, introduced in  \cite{FR}, while the  spaces $G\Gamma(p,\infty,w_1,w_2)$ appeared in \cite{GAG}.

If we assume that $uw_1(u)$ and $w_2$ are slowly varying functions, we can identify the Generalized Gamma space as an ${\mathcal L}$-space in the following way
$$G\Gamma(p,q,w_1,w_2)=L^{\mathcal L}_{p,(uw_1(u))^{1/q},L_q,(w_2(u))^{1/p},L_p} .$$
Thus, we can  apply the results from \S \ref{sereiteration} to interpolate two Generalized Gamma spaces with double weights.

\begin{cor}
Let $E$ be an r.i. space on $(0,1)$ and $\b, uw_1(u), w_2, uw_3(u), w_4\in SV(0,1)$. Let $1<p_0<p_1<\infty$, $1\leq q\leq \infty$ and consider the function 
$$\rho(u)=u^{\frac1{p_0}-\frac1{p_1}}\frac{w_2^{\frac1{p_0}}(u)\|(tw_1(t))^{\frac1{q_0}}\|_{\widetilde{L}_{q_0}(u,1)}}{w_4^{\frac1{p_1}}(u)\|(tw_3(t))^{\frac1{q_1}}\|_{\widetilde{L}_{q_1}(u,1)}},\quad u\in(0,1).$$

\begin{itemize}
\item[a)] If $0<\theta< 1$, then
$$\big(G\Gamma(p_0,q_0,w_1,w_2),G\Gamma(p_1,q_1,w_3,w_4)\big)_{\theta,\b,E}=L_{p,\B_\theta,E}$$
where  $\frac1{p}=\frac{1-\theta}{p_0}+\frac{\theta}{p_1}$ and $B_\theta(u)$ is equal to  $$\Big((w_2(u))^{\frac{1}{p_0}}\|(tw_1(t))^\frac{1}{q_0}\|_{\widetilde{L}_{q_0}(u,1)}\Big)^{1-\theta}\Big((w_4(u))^{\frac{1}{p_1}}\|(tw_3(t))^\frac{1}{q_1}\|_{\widetilde{L}_{q_1}(u,1)}\Big)^{\theta}\b(\rho(u)),$$
for $u\in(0,1)$.

\vspace{1mm}
\item[b)] If $\theta=0$ then
\begin{align*}
\big(G\Gamma(p_0,q_0,w_1,w_2),&G\Gamma(p_1,q_1,w_3,w_4)\big)_{0,\b,E}\\&=
L^{\mathcal L}_{p_0,B_0,E,(w_2(u))^{1/p_0},L_{p_0}}\cap L^{\mathcal L,\mathcal L}_{p_0,\b\circ\rho,E,(uw_1(u))^{1/q_0},L_{q_0},(w_2(u))^{1/p_0},L_{p_0}}
\end{align*}
where $B_0(u)=\|(tw_1(t))^\frac{1}{q_0}\|_{\widetilde{L}_{q_0}(u,1)}\b(\rho(u))$, $u\in(0,1)$.

\vspace{1mm}
\item[c)] If $\theta=1$ and  $\|\b\|_{\widetilde{E}(0,1)}<\infty$, then
$$
\big(G\Gamma(p_0,q_0,w_1,w_2),G\Gamma(p_1,q_1,w_3,w_4)\big)_{1,\b,E}=L^{\mathcal L,\mathcal R}_{p_1,\b\circ\rho,E,(uw_3(u))^{1/q_1},L_{q_1},(w_4(u))^{1/p_1},L_{p_1}}.
$$
\end{itemize}
\end{cor}

\subsection{$A$ and $B$-type spaces.}\hspace{1mm}

\vspace{2mm}
Finally, we consider the $A$ and $B$-type spaces studied by Pustylnik in \cite{pu2}.

\begin{defn}Given $1 <p<\infty$, $\alpha<1$ and $E$ an r.i. space on $(0,1)$. The space $A_{p,\alpha,E}$ is  the set of all measurable functions $f$ on $(0,1)$ such that
$$ \|f\|_{A_{p,\alpha,E}}=\Big\|\ell^{\alpha-1}(t)\int_t^1s^{\frac{1}{p}}f^{**}(s)\,\frac{ds}s\Big\|_
{\widetilde{E}(0,1)}<\infty$$
assumed that the function $(1+u)^{\alpha-1}$ belongs to $E$ (i.e. $\|\ell^{\alpha-1}(t)\|_{\widetilde{E}(0,1)}<\infty$ see \cite{pu2}).
The space $B_{p,\alpha,E}$ is the set of all measurable functions $f$ on $(0,1)$ such that 
$$
\|f\|_{B_{p,\alpha,E}}=\Big\|\sup_{0<s<t}s^{\frac{1}{p}}\ell^{\alpha-1}(s) f^{**}(s)\Big\|_{\widetilde{E}(0,1)}<\infty.$$
\end{defn}

The spaces of $B$-type when $\alpha=0$  first appeared in \cite{CwP}. General versions of these spaces were studied in \cite{PS1}. The main feature of the $A$ and $B$-type spaces is their optimality for  weak interpolation \cite{pu2,PS2}.

The $A$ and $B$-type spaces can be seen as $\mathcal{R}$ and $\mathcal{L}$-spaces, respectively. Indeed,
$$A_{p,\alpha,E}=L^{\mathcal{R}}_{p,\ell^{\alpha-1}(t),E,1,L_1}\quad B_{p,\alpha,E}=L^{\mathcal{L}}_{p,1,E,\ell^{\alpha-1}(t),L_\infty}.$$
Then, we can apply the results from \S \ref{sereiteration} to obtain the following interpolation formulae.

\vspace{1mm}

\begin{cor}
Let $E,\ E_0,\ E_1$ be r.i. spaces on $(0,1)$, $\b\in SV(0,1)$, $1< p_0<p_1< \infty$, $\alpha, \beta<1$ and assume that $(1+u)^{\alpha-1}$, $(1+u)^{\beta-1}$, belongs to $E_0$, $E_1$, respectively. Consider the function 
$$\rho(u)=u^{\frac{1}{p_0}-\frac{1}{p_1}}\frac{\|\ell^{\alpha-1}(t)\|_{\widetilde{E}_0(0,u)}}{\|\ell^{\beta-1}(t)\|_{\widetilde{E}_1(0,u)}},\qquad u \in (0,1).$$
\begin{itemize}

\vspace{1mm}
\item[a)] If $0<\theta< 1$, then
\begin{equation*}
\big(A_{p_0,\alpha,E_0},A_{p_1,\beta,E_1}\big)_{\theta,\b,E}=L_{p,\B_\theta,E},
\end{equation*}
where $\frac1{p}=\frac{1-\theta}{p_0}+\frac{\theta}{p_1}$ and $$\B_\theta(u)=\|\ell^{\alpha-1}(t)\|^{1-\theta}_{\widetilde{E}_0(0,u)}\|\ell^{\beta-1}(t)\|^{\theta}_{\widetilde{E}_1(0,u)}\b(\rho(u)),\quad u \in(0,1).$$

\vspace{1mm}
\item[b)] If $\theta=0$, then
$$
\big(A_{p_0,\alpha,E_0},A_{p_1,\beta,E_1}\big)_{0,\b,E}=L^{\mathcal R,\mathcal L}_{p_0,\b\circ\rho,E,\ell^{\alpha-1}(u),E_0,1,L_1}.$$

\vspace{1mm}
\item[c)] If $\theta=1$ and  $\|\b\|_{\widetilde{E}(0,1)}<\infty$, then
$$
\big(A_{p_0,\alpha,E_0},A_{p_1,\beta,E_1}\big)_{1,\b,E}=L^{\mathcal R}_{p_1,\B_1,E,1,L_{1}}\cap L^{\mathcal R,\mathcal R}_{p_1,\b\circ\rho,E,\ell^{\beta-1}(u),E_1,1,L_1}$$
where $\B_1(u)=\|\ell^{\beta-1}(t)\|_{\widetilde{E}_1(0,u)}\b(\rho(u))$, $u \in(0,1)$.
\end{itemize}
\end{cor}

\begin{cor}
Let $E,\, E_0,\, E_1$ be r.i. spaces on $(0,1)$, $\b\in SV(0,1)$, $1<p_0<p_1<\infty$ and $\alpha,\beta<1$. Consider the function $\rho(u)=u^{\frac{1}{p_0}-\frac{1}{p_1}}\frac{\ell^{\alpha-1}(u)\varphi_{E_0}(\ell(u))}{\ell^{\beta-1}(u)\varphi_{E_1}(\ell(u))}$, $u \in (0,1)$.

\vspace{1mm}
\begin{itemize}
\item[a)] If $0<\theta< 1$, then
$$\big(B_{p_0,\alpha,E_0},B_{p_1,\beta,E_1}\big)_{\theta,\b,E}=L_{p,\B_\theta,E},$$
where $\frac1{p}=\frac{1-\theta}{p_0}+\frac{\theta}{p_1}$ and $$\B_\theta(u)=\big(\ell^{\alpha-1}(u)\varphi_{E_0}(\ell(u))\big)^{1-\theta}\big(\ell^{\beta-1}(u)\varphi_{E_1}(\ell(u))\big)^{\theta}\b(\rho(u)),\quad u \in(0,1).$$

\vspace{1mm}
\item[b)] If $\theta=0$, then
$$
\big(B_{p_0,\alpha,E_0},B_{p_1,\beta,E_1}\big)_{0,\b,E}=L^{\mathcal L}_{p_0,\B_0,E,\ell^{\alpha-1}(u),L_\infty}\cap L^{\mathcal L,\mathcal L}_{p_0,\b\circ\rho,E,1,E_0,\ell^{\alpha-1}(u),L_\infty}
$$
where $B_0(u)=\varphi_{E_0}(\ell(u))\b(\rho(u))$, $u\in(0,1)$.

\vspace{1mm}
\item[c)] If $\theta=1$ and  $\|\b\|_{\widetilde{E}(0,1)}<\infty$, then
$$
\big(B_{p_0,\alpha,E_0},B_{p_1,\beta,E_1}\big)_{1,\b,E}=L^{\mathcal L,\mathcal R}_{p_1,\b\circ\rho,E,1,E_1,\ell^{\beta-1}(u),L_\infty}.$$
\end{itemize}
\end{cor}

\begin{cor}
Let $E,\, E_0,\, E_1$ be r.i. spaces on $(0,1)$, $\b\in SV(0,1)$, $1<p_0<p_1<\infty$ and $\alpha,\beta<1$. Consider the function $\rho(u)=u^{\frac{1}{p_0}-\frac{1}{p_1}}\frac{\|\ell^{\alpha-1}(t)\|_{\widetilde{E}_0(0,u)}}{\ell^{\beta-1}(u)\varphi_{E_1}(\ell(u))}$, $u \in (0,1)$.

\vspace{1mm}
\begin{itemize}
\item[a)] If $0<\theta< 1$, then
$$\big(A_{p_0,\alpha,E_0},B_{p_1,\beta,E_1}\big)_{\theta,\b,E}=L_{p,\B_\theta,E},$$
where $\frac1{p}=\frac{1-\theta}{p_0}+\frac{\theta}{p_1}$ and $$\B_\theta(u)=\|\ell^{\alpha-1}(t)\|^{1-\theta}_{\widetilde{E}_0(0,u)}\big(\ell^{\beta-1}(u)\varphi_{E_1}(\ell(u))\big)^{\theta}\b(\rho(u)),\quad u \in(0,1).$$

\vspace{1mm}
\item[b)] If $\theta=0$, then
$$
\big(A_{p_0,\alpha,E_0},B_{p_1,\beta,E_1}\big)_{0,\b,E}
=L^{\mathcal R,\mathcal L}_{p_0,\b\circ\rho,E,\ell^{\alpha-1}(u),E_0,1,L_1}.
$$

\vspace{1mm}
\item[c)] If $\theta=1$ and  $\|\b\|_{\widetilde{E}(0,1)}<\infty$, then
$$
\big(A_{p_0,\alpha,E_0},B_{p_1,\beta,E_1}\big)_{1,\b,E}=L^{\mathcal L,\mathcal R}_{p_1,\b\circ\rho,E,1,E_1,\ell^{\beta-1}(u),L_\infty}.$$
\end{itemize}
\end{cor}

Our last result completes \cite[Corollary 5.19]{FMS-RL2}.
\begin{cor}\label{corAB}
Let $E,\, E_0,\, E_1$ be r.i. spaces on $(0,1)$, $\b\in SV(0,1)$, $1<p_0<p_1<\infty$ and $\alpha,\beta<1$. Consider the function $
\rho(u)=u^{\frac{1}{p_0}-\frac{1}{p_1}}\frac{\ell^{\alpha-1}(u)\varphi_{E_0}(\ell(u))}{\|\ell^{\beta-1}(t)\|_{\widetilde{E}_1(0,u)}},\quad u \in (0,1).
$

\vspace{1mm}
\begin{itemize}
\item[a)] If $0<\theta< 1$, then
$$\big(B_{p_0,\alpha,E_0},A_{p_1,\beta,E_1}\big)_{\theta,\b,E}=L_{p,\B_\theta,E},$$
where $\frac1{p}=\frac{1-\theta}{p_0}+\frac{\theta}{p_1}$ and $$\B_\theta(u)=\big(\ell^{\alpha-1}(u)\varphi_{E_0}(\ell(u))\big)^{1-\theta}\|\ell^{\beta-1}(t)\|^\theta_{\widetilde{E}_0(0,u)}\b(\rho(u)),\quad u \in(0,1).$$

\vspace{1mm}
\item[b)] If $\theta=0$, then
\begin{equation}\label{einter5}
\big(B_{p_0,\alpha,E_0},A_{p_1,\beta,E_1}\big)_{0,\b,E}=L^{\mathcal L}_{p_0,\B_0,E,\ell^{\alpha-1}(u),L_\infty}\cap L^{\mathcal L,\mathcal L}_{p_0,\b\circ\rho,E,1,E_0,\ell^{\alpha-1}(u),L_\infty}
\end{equation}
where $B_0(u)=\varphi_{E_0}(\ell(u))\b(\rho(u))$, $u\in(0,1)$.

\vspace{1mm}
\item[c)] If $\theta=1$ and  $\|\b\|_{\widetilde{E}(0,1)}<\infty$, then
$$
\big(B_{p_0,\alpha,E_0},A_{p_1,\beta,E_1}\big)_{1,\b,E}=L^{\mathcal R}_{p_1,\B_1,E,1,L_1}\cap L^{\mathcal R,\mathcal R}_{p_1,\b\circ\rho,E,\ell^{\beta-1}(u),E_1,1,L_1},$$
where $B_1(u)=\|\ell^{\beta-1}(t)\|_{\widetilde{E}_1(0,u)}\b(\rho(u))$, $u\in(0,1)$.
\end{itemize}
\end{cor}

\begin{rem}
Comparing \eqref{einter5} with Theorem 5.19 b) from \cite{FMS-RL2} one can deduce that
\begin{align*}
L^{\mathcal L}_{p_0,B_{0},E,\ell^{\alpha-1}(u),L_\infty}\cap\big(L_{p_0,\ell^{\alpha-1}(u),L_\infty},&A_{p_1,\beta,E_1}\big)^{\mathcal L}_{0,\b(u \varphi_{E_{0}}(\ell(u))),E,1,E_0}\\  &=L^{\mathcal L}_{p_0,\B_0,E,\ell^{\alpha-1}(u),L_\infty}\cap L^{\mathcal L,\mathcal L}_{p_0,\b\circ\rho,E,1,E_0,\ell^{\alpha-1}(u),L_\infty}.
\end{align*}
\end{rem}

{\small \hspace{-5mm}{\textbf{Acknowledgments.}}
The authors have been partially supported by grant MTM2017-84058-P (AEI/FEDER, UE).
The second author also thanks the Isaac Newton Institute for Mathematical Sciences, Cambridge, for support and hospitality during the programme \emph{Approximation, Sampling and Compression in Data Science} where the work on this paper was undertaken; this work was supported by EPSRC grant no. EP/R014604/1. Finally, the second author thanks Óscar Domínguez for useful conversations at the early stages of this work, and for pointing out the reference \cite{FFGKR}. }

\end{document}